\documentclass[12pt]{article}
\usepackage{amsmath}
\usepackage{amssymb}
\usepackage{graphicx}
\usepackage[numbers,sort&compress]{natbib}
\usepackage[colorlinks,
            anchorcolor=blue,
            citecolor=green]
            {hyperref}
\title{Well-posedness to nonlinear Schr\"odinger-Gerdjikov-Ivanon equation}

\author{Sucai Niu, 
Junyi Zhu\thanks{E-mail: jyzhu@zzu.edu.cn}\\
\scriptsize{\sl School of Mathematics and Statistics, Zhengzhou University, Zhengzhou, Henan 450001, China}\\
}

\date{}
\begin{document}
\newtheorem{theorem}{Theorem}
\newtheorem{definition}{Definition}
\newtheorem{lemma}{Lemma}
\newtheorem{corollary}{Corollary}
\newtheorem{proposition}{Proposition}
\newtheorem{remark}{Remark}
\hypersetup{CJKbookmarks=true}

\maketitle
\begin{abstract}
The Riemann-Hilbert approach is extended to discuss the well-posedness of the nonlinear Schr\"odinger-Gerdjikov-Ivanon equation. The Lipschitz continuity of potential in $H^{2}(\mathbb{R})\cap H^{1,1}(\mathbb{R})$ to scattering data is obtained through direct scattering transform. Two Riemann-Hilbert problems are constructed, and two sets of the reflection coefficients, that is $r(k)$ and $r_\pm(z)$, are introduced. The Lipschitz continuity from the reflection coefficients $r_\pm(z)$ in $H^{1}(\mathbb{R})\cap L^{2,1}(\mathbb{R})$ to the potential is estimated via the potential reconstruction. Existence of global solutions of NLS-GI equation is considered by the Riemann-Hilbert problem without eigenvalues or resonances.


\end{abstract}

\section{Introduction}

The nonlinear Schr\"odinger equation (NLS) and the Gerdjikov-Ivanov equation (GI) are important integrable models, with profound implications in mathematical physics. The NLS equation, known for its simple cubic nonlinear structure and rich soliton solutions, serves as a paradigm for studying phenomena such as wave packet self-focusing and modulation instability, has important applications in nonlinear optics \cite{AG1989,chaos10-471}, deep water wave theory \cite{sam46-133}, Bose-Einstein condensates \cite{jpa15-2599}. The GI equation, as one of important derivative NLS equation, has vital applications in different fields, such as long-wavelength dynamics of dispersive Alfv\'en waves in the plasma physics field \cite{pf14-2733}, the subpicosecond or femtosecond pulse in a single-mode fiber \cite{ol12-516}, the propagation of nonlinear pulses in optical fibers \cite{AG2001}, and so on.

The combined NLS-GI equation
\begin{equation}\label{a1}
u_{t}=iu_{xx}-2i|u|^{2}u+u^{2}\overline{u}_{x}+\frac{i}{2}|u|^{4}u,\\
\end{equation}
introduces more complex nonlinear terms and dispersion relations, revealing the influence of higher-order effects on soliton dynamics. This equation and its nonlocal forms \cite{jmp65-103501} has been well studied \cite{mplb32-1850088,tmp209-1537,ps98-125267,eajam15-163}.
The well-posedness for derivative NLS equation with the inverse scattering transform method have been studied in \cite{cpde41-271,imrn2018-5663,dpde14-271,cpde43-1151,jpam78-33}, and other integrable equations \cite{jde366-320,sam152-111,camb45-497,jde414-34,jde440-113446}.
In this paper, we extend the inverse scattering transform to discuss the well-posedness for the NLS-GI equation \eqref{a1}. Here, we discuss the well-posedness of the Cauchy problem for the NLS-GI equations \eqref{a1}
with initial condition $u(x,0)=u_{0}(x)\in H^2(\mathbb{R})\cap H^{1,1}(\mathbb{R})$.

The NLS-GI equation is equivalent to the compatibility condition of the following linear spectral problem
\begin{equation}\label{a2}
Y_{x}=UY, \quad
U=\left( {{\begin{array}{*{20}c}
i(-\lambda+\frac{1}{2}|u|^2)& i(1-\lambda)u  \\
-i\overline{u} &  i(\lambda-\frac{1}{2}|u|^2)\\
\end{array} }} \right),
\end{equation}
and
\begin{equation}\label{a3}
\begin{array}{l}
Y_{t}=VY,\quad V=\left( {{\begin{array}{*{20}c}
A & (1-\lambda)(2i\lambda u-u_{x}) \\
-2i\lambda\overline{u}-\overline{u}_{x} & -A \\
\end{array} }} \right),
\end{array}
\end{equation}
where
$$A=-2i\lambda^{2}+i\lambda |u|^2+\frac{1}{2}(u\overline{u}_{x}-u_{x}\overline{u})+\frac{i}{4}|u|^{4}-i|u|^2.$$
We note that the inverse scattering approach to the combine NLS-GI \eqref{a1} is more challenging due to the spectral parameter $\lambda=z+1=k^2+1$ be considered.

The article is organized as follows. In section 2, the analytical and asymptotic properties of the Jost are discussed by introducing two transformation matrices. The Lipschitz continuity from the potential in $H^{2}(\mathbb{R})\cap H^{1,1}(\mathbb{R})$ to scattering data is discussed. In section 3, the Riemann-Hilbert problems and their solutions equivalent to spectral problem \eqref{a1} are presented, and prior estimates through scattering data for $x\in\mathbb{R}^{+}$ and $x\in\mathbb{R}^{-}$ are given. In section 4, potential reconstruction is obtained by virtue of the Riemann-Hilbert problem, and well-posedness analysis and the Lipschitz continuity of the maps from reflection coefficients to the potential are presented via the theories of Fourier transform in detail. In the section 5, the existence of global solutions is shown.

\setcounter{equation}{0}
\section{The direct scattering transform and Lipschitz continuity}\label{sec2}
For convenience, we let $z=k^2\in\mathbb{R}, \lambda=z+1$ and introduce two modified Jost functions
\begin{equation}\label{mtb4}
 m_\pm(x;z)=T_1{Y}(x;z)\mathrm{e}^{i(z+1)x}, \quad n_\pm(x;z)=T_2{Y}(x;z)\mathrm{e}^{-i(z+1)x},
\end{equation}
with the boundary conditions
\begin{equation}\label{mtb5}
 m_\pm(x;z)\to e_1=\left(\begin{array}{c}
 1\\
 0
 \end{array}\right), \quad n_\pm(x;z)\to e_2=\left(\begin{array}{c}
 0\\
 1
 \end{array}\right), \quad x\to\pm\infty,
\end{equation}
where
\begin{equation}\label{mtb1}
 T_1=\left(\begin{matrix}
 1&0\\
 i\overline{u}(x)&-2ik^2
 \end{matrix}\right), \quad T_2=\left(\begin{matrix}
 -2k&-ku(x)\\
 0&-ik
  \end{matrix}\right),
\end{equation}
Then the Jost functions $m_\pm(x;z)$ and $n_\pm(x;z)$ satisfy the Volterra integral equations
\begin{equation}\label{mtb6}
 m_\pm(x;z)=e_1+\int_{\pm\infty}^x\left(\begin{matrix}
 1&0\\
 0&\mathrm{e}^{2i(z+1)(x-y)}
 \end{matrix}\right)Q_1(y)m_\pm(y;z)\mathrm{d}y,
\end{equation}
and
\begin{equation}\label{mtb7}
 n_\pm(x;z)=e_2+\int_{\pm\infty}^x\left(\begin{matrix}
 \mathrm{e}^{-2i(z+1)(x-y)}&0\\
 0&1
 \end{matrix}\right)Q_2(y)n_\pm(y;z)\mathrm{d}y,
\end{equation}
where
\[Q_1(x)=\left(\begin{matrix}
 0&\frac{u}{2}\\
 \bar{w}&0
 \end{matrix}\right),  \quad
 Q_2(x)=\left(\begin{matrix}
 0& w\\
 \frac{\overline{u}}{2}&0
 \end{matrix}\right), \quad w=-iu_x+2u-\frac{1}{2}|u|^2u.\]

Since ${\rm Im}(z+1)={\rm Im}z, z\in\mathbb{C}$, then for every $x\in\mathbb{R}$, via the Neumann series, we can show that $m_\mp(x;\cdot)$ and $n_\pm(x;\cdot)$ are continued analytically in $\mathbb{C}^\pm=\{z:{\rm Im}z\gtrless0\}$, if $u\in L^3(\mathbb{R})\cap L^1(\mathbb{R}), u_x\in L^1(\mathbb{R})$. In addition, for every $z\in\mathbb{R}$, $m_\pm$ and $n_\pm$ are bounded in the space $\sup\limits_{z\in\mathbb{R}}\|\cdot\|_{L_x^\infty}$.

According to the Volterra integral equations \eqref{mtb6} and \eqref{mtb7}, we introduce two integral operators
\begin{equation}\label{nlsgia1}
({\cal{K}}_{1\pm}f)(x;z)=\int_{\pm\infty}^x\left(\begin{matrix}
 1&0\\
 0&\mathrm{e}^{2i(z+1)(x-y)}
 \end{matrix}\right)W_1(y)f(y;z)\mathrm{d}y,
\end{equation}
and
\begin{equation}\label{nlsgia2}
({\cal{K}}_{2\pm}f)(x;z)=\int_{\pm\infty}^x\left(\begin{matrix}
 \mathrm{e}^{-2i(z+1)(x-y)}&0\\
 0&1
 \end{matrix}\right)W_2(y)f(y;z)\mathrm{d}y,
\end{equation}
for every vector function $f(x;z)\in L_x^\infty(\mathbb{R},L_z^2(\mathbb{R}))$. Here and after, for every matrix or vector $A$, we take
\[\|A\|_{L_z^p}=\||A|\|_{L_z^p}, \quad |A|=\sqrt{\mathrm{tr}(A^\dagger A)}.\]
Using the Minkowski inequality in integral form, we find that
\begin{equation}\label{nlsgia3}
 \begin{aligned}
   \|({\cal{K}}_{j\pm}f)\|\|_{L_z^2}&\leq\mp\int_{\pm\infty}^x|W_j(y)|\|f(y,z)\|_{L_z^2}\mathrm{d}y\\
  &\leq\|f(x,z)\|_{L_x^\infty L_z^2}\left(\mp\int_{\pm\infty}^x|W_j(y)|\mathrm{d}y\right), \quad j=1,2,
 \end{aligned}
\end{equation}
which imply that
\begin{equation}\label{nlsgia4}
 \|({\cal{K}}_{j\pm}^nf)\|\|_{L_z^2}\leq\frac{\|W_j\|_{L^1}^n}{n!}\|f(x;z)\|_{L_x^\infty L_z^2}, \quad j=1,2,
\end{equation}
in terms of iteration.  We note that, for $j=1,2$,
\begin{equation}\label{nlsgia4a}
\begin{aligned}
\|W\|_{L^p}:&=\|W_j\|_{L^p}\lesssim \|w\|_{L^p}+\|u\|_{L^p}, \quad p\geq1, \\ \|w\|_{L^1}&\lesssim\|u_x\|_{L^1}+\|u\|_{L^1}+\|u\|_{L^3}^3,
\end{aligned}
\end{equation}
if $u\in L^3(\mathbb{R})\cap L^1(\mathbb{R}), u_x\in L^1(\mathbb{R})$. In addition, if $u\in H^1$, we also have
\begin{equation}\label{nlsgia4b}
\begin{aligned}
\|w\|_{L^2}&\lesssim\|u_x\|_{L^2}+\|u\|_{L^2}+\|u\|_{L^6}^3\\
&\lesssim \|u_x\|_{L^2}+\|u\|_{L^2}+\|u\|_{L^2}^2\|u_x\|_{L^2},
\end{aligned}
\end{equation}
in view of the Sharp Gagliardo-Nirenberg (GN) inequality \cite{cmp87-567,jdde18-1069}
\begin{equation}\label{mtb27b}
\|u\|_{L^p(\mathbb{R})}\lesssim\|u\|^{1-\theta}_{L^2(\mathbb{R})}\|u_x\|^{\theta}_{L^2(\mathbb{R})}, \quad \theta=\frac{p-2}{2p}.
\end{equation}

From \eqref{nlsgia4}, we can get that the operators $I-{\cal{K}}_{j\pm}$ is invertible in the space $L_x^\infty(\mathbb{R},L_z^2(\mathbb{R}))$, and
\begin{equation}\label{nlsgia4c}
\|(I-{\cal{K}}_{j\pm})^{-1}\|_{L_x^\infty L_z^2\mapsto L_x^\infty L_z^2}\leq \mathrm{e}^{\|W_j\|_{L^1}}=\mathrm{e}^{\|W\|_{L^1}}.
\end{equation}

Using the integral operators ${\cal{K}}_{j\pm}$, we rewrite the Volterra integral equation \eqref{nlsgia1} and \eqref{nlsgia2} as the contract form
\begin{equation}\label{nlsgia5}
m_\pm(x;z)=e_1+{\cal{K}}_{1\pm}m_\pm(x;z), \quad n_\pm(x;z)=e_2+{\cal{K}}_{2\pm}n_\pm(x;z).
\end{equation}
Furthermore, we can get
\begin{equation}\label{nlsgia6}
m_\pm(x;z)-e_1=(I-{\cal{K}}_{1\pm})^{-1}{\cal{K}}_{1\pm}e_1,
\end{equation}
and
\begin{equation}\label{nlsgia7}
n_\pm(x;z)-e_2=(I-{\cal{K}}_{2\pm})^{-1}{\cal{K}}_{2\pm}e_2.
\end{equation}

We note that
\begin{equation}\label{nlsgia8}
\|{\cal{K}}_{j\pm}e_j\|_{L_z^2}=\sqrt{\pi}\left(\mp\int_{\pm\infty}^x|w(y)|^2dy\right)^{\frac{1}{2}}, \quad j=1,2,
\end{equation}
in terms of the Plancherel formula. In addition, we have
\begin{equation}\label{nlsgia9}
 \|{\cal{K}}_{j\pm}e_j\|_{L_x^\infty L_z^2}\lesssim\|w\|_{L^2}, \quad j=1,2.
\end{equation}
\begin{lemma}\label{lema1}
If $u\in H^{1,1}$ , then $m_\pm(x;z)-e_1, n_\pm(x;z)-e_2\in L_x^\infty(\mathbb{R}^\pm,H_z^1(\mathbb{R}))$.
Moreover, if $u\in H^{1,1}$, then $m_\pm(x;z)-e_1, n_\pm(x;z)-e_2\in L_x^2(\mathbb{R}^\pm,L_z^2(\mathbb{R}))$.
\end{lemma}
{\bf Proof}. We give the proof for $m_+(x;z)$. From \eqref{nlsgia4}, \eqref{nlsgia9} and \eqref{nlsgia6}, we find that
\begin{equation}\label{nlsgia10}
 \|m_+(x;z)-e_1\|_{L_x^\infty L_z^2}\leq \mathrm{e}^{\|W\|_{L^1}}\|w\|_{L^2},
\end{equation}
where $\|W\|_{L^1}$ is define in \eqref{nlsgia4a}. Since $|u|_{L^1}\leq\sqrt{\pi}|u|_{L^{2,1}}$, then from \eqref{nlsgia4a}, \eqref{nlsgia4b} and the sharp GN inequality \eqref{mtb27b}, we can get that if $u\in H^{1,1}$ then $m_+(x;z)-e_1\in L_x^\infty(\mathbb{R}^+,L_z^2(\mathbb{R}))$.

Next, we prove $m_+(x;z)-e_1\in L_x^2(\mathbb{R}^+,L_z^2(\mathbb{R}))$. we note that equation \eqref{nlsgia8} implies, for every $x\in \mathbb{R}^+$, that
\begin{equation}\label{nlsgia11}
\begin{aligned}
\|{\cal{K}}_{1+}e_1\|_{L_x^2 L_z^2}&\lesssim\left(\int_0^\infty\int_x^\infty|w(y)|^2\mathrm{d}y\mathrm{d}x\right)^{\frac{1}{2}}\\
&\leq\left(\int_0^\infty|x||w(x)|^2\mathrm{d}x\right)^{\frac{1}{2}}
\leq\|w\|_{L^{2,1}}^{\frac{1}{2}}\|w\|_{L^2}^{\frac{1}{2}},
\end{aligned}
\end{equation}
via integration by parts. In addition, from \eqref{nlsgia1}, we can get
\begin{equation}\label{nlsgia12}
 \|{\cal{K}}_{1+}(m_+(x;z)-e_1)\|_{L_z^2}\leq \|m_+(x;z)-e_1\|_{L_x^\infty L_z^2}\int_x^\infty|W(y)|\mathrm{d}y.
\end{equation}
Then for every $x\in\mathbb{R}^+$, from \eqref{nlsgia10} and \eqref{nlsgia12},  we have
\begin{equation}\label{nlsgia13}
  \begin{aligned}
  &\|{\cal{K}}_{1+}(m_+(x;z)-e_1)\|^2_{L_x^2 L_z^2}\\
  &\lesssim \int_0^\infty\left(\int_x^\infty|W(y)|\mathrm{d}y\right)^2\mathrm{d}x
 =\left(\int_0^1+\int_1^\infty\right)\left(\int_x^\infty|W(y)|\mathrm{d}y\right)^2\mathrm{d}x\\
 &\quad=\int_0^1\left(\big(\int_0^\infty-\int_0^x\big)|W(y)|\mathrm{d}y\right)^2\mathrm{d}x
 +\int_1^\infty\left(\int_x^\infty|W(y)|\mathrm{d}y\right)^2\mathrm{d}x\\
 &\leq\int_0^1\left(\int_0^\infty|W(y)|\mathrm{d}y\right)^2\mathrm{d}x+\int_0^1\left(\int_0^x|W(y)|\mathrm{d}y\right)^2\mathrm{d}x\\
 &\quad +\int_1^\infty\left(\int_x^\infty|W(y)|\mathrm{d}y\right)^2\mathrm{d}x\\
 &\leq\|W\|_{L^1}^2+\left(\int_0^1|W(x)|dx\right)^2+\|W\|_{L^{2,1}}^2,
  \end{aligned}
\end{equation}
where we have used integration by parts and the inequality (\cite{HLP1988}, Theorem 328)
\[\int_1^\infty\left(\int_x^\infty|W(y)|\mathrm{d}y\right)^2\mathrm{d}x\leq\int_1^\infty|xW(x)|^2dx.\]

Since equation \eqref{nlsgia5} gives
\begin{equation}\label{nlsgia14}
\|m_+(x;z)-e_1\|_{L_x^2 L_z^2}\leq \|{\cal{K}}_{1+}e_1\|_{L_x^2 L_z^2}+\|{\cal{K}}_{1+}(m_+(x;z)-e_1)\|_{L_x^2 L_z^2}.
\end{equation}
Hence, if $u\in H^{1,1}$, then from \eqref{nlsgia11}, \eqref{nlsgia13} and \eqref{nlsgia14}, we prove that $m_+(x;z)-e_1\in L_x^2(\mathbb{R}^+,L_z^2(\mathbb{R}))$.

To finish the proof of the first statement. we need to prove $\partial_zm_+(x,z)\in L_x^\infty(\mathbb{R}^+,L_z^2(\mathbb{R}))$. To this end, differentiating both sides of the equation \eqref{nlsgia5} of $m_+$  with respect to $z$, we can obtain that
\begin{equation}\label{nlsgia15}
(I-{\cal{K}}_{1+})\partial_zm_+(x;z)=\partial_z({\cal{K}}_{1+}e_1)+(\partial_z{\cal{K}}_{1+})(m_+(x;z)-e_1).
\end{equation}
It is easy to see, for $0<x<y$, that
\begin{equation}\label{nlsgia16}
\|\partial_z({\cal{K}}_{1+}e_1)\|_{L_x^\infty L_z^2}\leq 2\sqrt{\pi}\left(\int_x^\infty|yw(y)|^2\mathrm{d}y\right)^{\frac{1}{2}}\lesssim\|w\|_{L^{2,1}}.
\end{equation}
In addition, for $0<x<y$, we find that
\begin{equation}\label{nlsgia17}
\begin{aligned}
\|(\partial_z{\cal{K}}_{1+})(m_+(x;z)-e_1)\|_{L_z^2}&\leq\int_x^\infty|yw(y)|\|m_+(y;z)-e_1\|_{L_z^2}\mathrm{d}y\\
&\leq\|w\|_{L^{2,1}}\|m_+(x;z)-e_2\|_{L_x^2L_z^2}.
\end{aligned}
\end{equation}
Thus $\partial_zm_+(x,z)\in L_x^\infty(\mathbb{R}^+,L_z^2(\mathbb{R}))$ is proved from \eqref{nlsgia15}-\eqref{nlsgia17} via \eqref{nlsgia14} and the bound \eqref{nlsgia4} of the operator $(I-{\cal{K}}_{1+})^{-1}$.

Hence, statements for $m_+(x;z)$ is proved. The other cases can be proved similarly. \quad$\Box$

The asymptotic behaviors of the Jost functions $m_\pm(x;z)$ and $n_\pm(x;z)$ near $z\to\infty$ take the following form
\begin{equation}\label{mtb10}
\lim\limits_{z\to\infty}m_\pm(x;z)=e_1, \quad
\lim\limits_{z\to\infty}n_\pm(x;z)=e_2.
\end{equation}
and
\begin{equation}\label{mtb12}
\lim\limits_{z\to\infty}2iz[m_\pm(x;z)-e_1]=-q_{\pm}(x) e_1-\overline{w(x)} e_2,
\end{equation}
\begin{equation}\label{mtb13}
\lim\limits_{z\to\infty}2iz[n_\pm(x;z)-e_2]={w(x)} e_1+\overline{q_{\pm}(x)} e_{2},
\end{equation}
where
\begin{equation}\label{mtb14}
q_{\pm}(x)=\frac{1}{2}\int_{\pm\infty}^{x}u(y)\overline{w(y)}\mathrm{d}y, \quad w=-iu_x+2u-\frac{1}{2}|u|^2u.
\end{equation}
\begin{lemma}
If $u\in H^{1,1}(\mathbb{R})\cap H^2(\mathbb{R})$, then
\begin{equation}\label{mtb15}
 \begin{aligned}
 m_\pm^{\langle2\rangle}(x;z):=2iz[m_\pm(x;z)-e_1]+q_{\pm}(x) e_1+\overline{w(x)} e_2\in L_x^\infty(\mathbb{R}^\pm,L_z^2(\mathbb{R})),\\
n_\pm^{\langle2\rangle}(x;z):=2iz[n_\pm(x;z)-e_2]-{w(x)} e_1-\overline{q_{\pm}(x)} e_{2}\in L_x^\infty(\mathbb{R}^\pm,L_z^2(\mathbb{R})).
 \end{aligned}
\end{equation}
\end{lemma}
{\bf Proof}~~We give the proof only for the case of $m_+^{\langle2\rangle}(x;z)$. Using \eqref{nlsgia6}, we find that
\begin{equation}\label{mtb16}
\begin{aligned}
(I-{\cal{K}}_{1+})m_+^{\langle2\rangle}(x;z)=&2iz{\cal{K}}_{1+}e_1-{\cal{K}}_{1+}q_{\pm}(x) e_1+\overline{w(x)} e_2\\
&+q_{\pm}(x) e_1-{\cal{K}}_{1+}\overline{w(x)} e_2.
\end{aligned}
\end{equation}
It is noted that
\[\begin{aligned}
&q_1e_1-{\cal{K}}_{1+}\overline{w(x)} e_2=0,\\
&2i(z+1){\cal{K}}_{1+}e_1+\overline{w(x)} e_2=\int_{+\infty}^xe^{2i(z+1)(x-y)}\overline{w_y(y)}\mathrm{d}ye_2.
\end{aligned}\]
Substituting them into \eqref{mtb16}, we can get that
\[
\begin{aligned}
(I-{\cal{K}}_{1+})m_+^{\langle2\rangle}(x;z)
=\int_{+\infty}^xe^{2i(z+1)(x-y)}[\overline{w_y(y)}-2i\overline{w(y)}-\overline{w(y)}q_+(y)]\mathrm{d}ye_2.
\end{aligned}
\]
By virtue of the bound \eqref{nlsgia4} and the Plancherel formula, we have
\begin{equation}\label{mtb17}
\|m_+^{\langle2\rangle}(x;z)\|_{L_x^\infty L_z^2}\leq \mathrm{e}^{\|Q\|_{L^1}}(\|w_x\|_{L^2}+2\|w\|_{L^2}+\|q_+\|_{L^\infty}\|w\|_{L^2}),
\end{equation}
in view of $\|q_+\|_{L^\infty}\leq\frac{1}{2}\|u\|_{L^2}\|w\|_{L^2}$ and
\[\|w_x\|_{L^2}\leq\|u_{xx}\|_{L^2}+2\|u_x\|_{L^2}+2\|u\|_{L^\infty}^2\|u_x\|_{L^2}.\]
Since $\|u\|_{L^\infty}\leq\frac{1}{\sqrt{2}}\|u\|_{H^1}$, then we prove that $m_+^{\langle2\rangle}\in L_x^\infty(\mathbb{R}^\pm,L_z^2(\mathbb{R}))$ if $u\in H^{1,1}\cap H^2$.
The other cases can be proved similarly. \qquad $\Box$

\begin{lemma}\label{lem4}
The maps
\begin{equation}\label{mtb36}
 H^{1,1}(\mathbb{R})\ni u\mapsto[m_\pm(x;z)-e_1,n_\pm(x;z)-e_2]\in L^\infty(\mathbb{R}^\pm;H_z^1(\mathbb{R})),
\end{equation}
and
\begin{equation}\label{mtb36a}
 H^{1,1}(\mathbb{R})\cap H^2(\mathbb{R})\ni u\mapsto[{m}_\pm^{\langle2\rangle}(x;z),{n}_\pm^{\langle2\rangle}(x;z)]\in L^\infty(\mathbb{R}^\pm;L_z^2(\mathbb{R})),
\end{equation}
are Lipschitz continuous.
\end{lemma}
{\bf Proof}~~ We give the proof for the case about $m_+(x,z)$. To this end, we let $\|u\|_{H^{1,1}(\mathbb{R})},\|\tilde{u}\|_{H^{1,1}(\mathbb{R})}\leq B$, for some $B$, $m_+(x;z)$ and $\tilde{m}_+(x;z)$ are the associated Jost functios.
From \eqref{nlsgia6}, we consider
\begin{equation}\label{mtb37}
\begin{aligned}
(m_+ -\tilde{m}_+)&=(I-{\cal{K}}_{1+})^{-1}({\cal{K}}_{1+}-{\cal{\tilde{K}}}_{1+})e_1\\
&\quad+(I-{\cal{K}}_{1+})^{-1}({\cal{K}}_{1+}-{\cal{\tilde{K}}}_{1+})(I-{\cal{\tilde{K}}}_{1+})^{-1}{\cal{\tilde{K}}}_{1+}e_1.
\end{aligned}
\end{equation}
From \eqref{nlsgia4} and \eqref{nlsgia9}, we can get
\begin{equation}\label{nlsgia18}
\begin{aligned}
\|(I-{\cal{K}}_{1+})^{-1}({\cal{K}}_{1+}-{\cal{\tilde{K}}}_{1+})e_1\|_{L_x^\infty L_z^2}
&\lesssim \mathrm{e}^{\|W\|_{L^1}}\|({\cal{K}}_{1+}-{\cal{\tilde{K}}}_{1+})e_1\|_{L_x^\infty L_z^2}\\
&\leq C(B)\|w-\tilde{w}\|_{L^2},
\end{aligned}
\end{equation}
and
\begin{equation}\label{nlsgia19}
\begin{aligned}
&\|(I-{\cal{K}}_{1+})^{-1}({\cal{K}}_{1+}-{\cal{\tilde{K}}}_{1+})(I-{\cal{\tilde{K}}}_{1+})^{-1}{\cal{\tilde{K}}}_{1+}e_1\|_{L_x^\infty L_z^2}\\
&\lesssim \mathrm{e}^{\|W\|_{L^1}}\|({\cal{K}}_{1+}-{\cal{\tilde{K}}}_{1+})(I-{\cal{\tilde{K}}}_{1+})^{-1}{\cal{\tilde{K}}}_{1+}e_1\|_{L_x^\infty L_z^2}\\
&\leq C(B)\|Q_1-\tilde{Q}_1\|_{L^1},
\end{aligned}
\end{equation}
where $C(B)$ is a positive constant. Moreover, the bounds \eqref{nlsgia4a} and \eqref{nlsgia4b} imply that
\[\|w-\tilde{w}\|_{L^2}\leq C(B)\|u-\tilde{u}\|_{H^1},\quad \|Q_1-\tilde{Q}_1\|_{L^1}\leq C(B)\|u-\tilde{u}\|_{H^{1,1}}.\]

Thus we have
\begin{equation}\label{nlsgia20}
 \|m_+ -\tilde{m}_+\|_{L_x^\infty L_z^2}\leq C(B)\|u-\tilde{u}\|_{H^{1,1}}.
\end{equation}

Similarly, using \eqref{nlsgia15}, we can get
\begin{equation}\label{nlsgia21}
\begin{aligned}
&\partial_zm_+-\partial_z\tilde{m}_+\\
&=(I-{\cal{K}}_{1+})^{-1}(\partial_z{\cal{K}}_{1+}e_1-\partial_z{\cal{\tilde{K}}}_{1+}e_1)\\
&\quad+(I-{\cal{K}}_{1+})^{-1}(\partial_z{\cal{K}}_{1+})(m_+-\tilde{m}_+)\\
&\quad+(I-{\cal{K}}_{1+})^{-1}({\cal{K}}_{1+}-{\cal{\tilde{K}}}_{1+})(I-{\cal{\tilde{K}}}_{1+})^{-1}\partial_z{\cal{\tilde{K}}}_{1+}e_1\\
&\quad+(I-{\cal{K}}_{1+})^{-1}({\cal{K}}_{1+}-{\cal{\tilde{K}}}_{1+})(I-{\cal{\tilde{K}}}_{1+})^{-1}(\partial_z{\cal{\tilde{K}}}_{1+})(\tilde{m}_+-e_1).
\end{aligned}
\end{equation}
From \eqref{nlsgia16} and \eqref{nlsgia4}, we can get
\[\|(I-{\cal{K}}_{1+})^{-1}(\partial_z{\cal{K}}_{1+}e_1-\partial_z{\cal{\tilde{K}}}_{1+}e_1)\|_{L_x^\infty L_z^2}\leq C(B)|w-\tilde{w}\|_{L^{2,1}}.\]
Using \eqref{nlsgia15}-\eqref{nlsgia17}, and taking similar estimation to  \eqref{nlsgia21}, we can obtain
\begin{equation}\label{nlsgia20a}
 \|\partial_zm_+ -\partial_z\tilde{m}_+\|_{L_x^\infty L_z^2}\leq C(B)\|u-\tilde{u}\|_{H^{1,1}}.
\end{equation}

So the Lipshitz continous for $u$ to $m_+-e_1$ in \eqref{mtb36} is proved.

The Lipshitz continous for $u$ to ${m}_+^{\langle2\rangle}$ in \eqref{mtb36a} can be prove via \eqref{mtb17} in a similar way. The other case of Lipshitz continous can be proved similarly. \qquad $\square$

Next, we define new Jost functions $\varphi_\pm(x;k)$ and $\phi_\pm(x;k)$ with $m_\pm(x;z)$ and $n_\pm(x;z)$ as
\begin{equation}\label{mtb44a}
\begin{aligned}
\varphi^{(1)}_{\pm}(x;k)=m^{(1)}_{\pm}(x;z),\quad \varphi^{(2)}_{\pm}(x;k)=\frac{1}{2k}[-i\overline{u}(x)m^{(1)}_{\pm}(x;z)+m^{(2)}_{\pm}(x;z)], \\
\phi^{(1)}_{\pm}(x;k)=\frac{-1}{2k}[n^{(1)}_{\pm}(x;z)+iu(x)n^{(2)}_{\pm}(x;z)], \quad\phi^{(2)}_{\pm}(x;k)=n^{(2)}_{\pm}(x;z),
\end{aligned}
\end{equation}
with the boundary conditions
\begin{equation}\label{mtb44}
 \varphi_\pm(x;k)\to e_1, \quad \phi_\pm(x;k)\to e_2, \quad x\to\pm\infty.
\end{equation}
Then $\varphi_\pm(x;k)$ and $\phi_\pm(x;k)$ are solutions of the Volterra integral equations
\begin{equation}
\label{eq44}
\varphi_{\pm}(x;k)=e_{1}+\int^{x}_{\pm\infty}\left(\begin{matrix}
1   &  0\\
0   &   \mathrm{e}^{2i(k^{2}+1)(x-y)}
\end{matrix}\right)Q(y)\varphi_{\pm}(y;k)\mathrm{d}y,
\end{equation}
\begin{equation}
\label{eq45}
\phi_{\pm}(x;k)=e_{2}+\int^{x}_{\pm\infty}\left(\begin{matrix}
 \mathrm{e}^{-2i(k^{2}+1)(x-y)} &  0\\
0   &   1
\end{matrix}\right)Q(y)\phi_{\pm}(y;k)\mathrm{d}y.
\end{equation}
where
\[Q(x;k)=\left( {{\begin{array}{*{20}c}
\frac{i}{2}|u|^{2}& ku  \\
-k\overline{u} &  -\frac{i}{2}|u|^{2}\\
\end{array} }} \right).\]

It is noted that $\varphi_\pm(x;k)$ and $\phi_\pm(x;k)$ satisfy the symmetry condition
\begin{equation}\label{mtb47}
\phi_\pm(x;k)=\sigma_1\sigma_3\overline{\varphi_\pm(x;\bar{k})}, \quad \sigma_1=\left(\begin{matrix}
0&1\\
1&0
\end{matrix}\right), \quad \sigma_1=\left(\begin{matrix}
1&0\\
0&-1
\end{matrix}\right).
\end{equation}
In addition, the Jost functions $\varphi_\pm(x;k)$ and $\phi_\pm(x;k)$ are not independent, and satisfy the following relation
\begin{equation}\label{mtb45}
\begin{aligned}
\varphi_-(x;k)=a(k)\varphi_+(x;k)+b(k)\mathrm{e}^{2i(k^2+1)x}\phi_+(x;k),\\
\phi_-(x;k)=-\overline{b(\bar{k})}\mathrm{e}^{-2i(k^2+1)x}\varphi_+(x;k)+\overline{a(\bar{k})}\phi_+(x;k).
\end{aligned}
\end{equation}
Substituting the Volterra integral equations \eqref{eq44} and \eqref{eq45} into \eqref{mtb45}, we can find that
\begin{equation}\label{mtb43}
\begin{aligned}
a(k)&=1+\int_{\mathbb{R}}[\frac{i}{2}|u(y)|^2\varphi_-^{(1)}(y;k)+ku(y)\varphi_-^{(2)}(y;k)]\mathrm{d}y,\\
b(k)&=-\int_{\mathbb{R}}e^{-2i(k^2+1)y}[k\overline{u(y)}\varphi_-^{(1)}(y;k)-\frac{i}{2}|u(y)|^2\varphi_-^{(2)}(y;k)]\mathrm{d}y,
\end{aligned}
\end{equation}
which imply that $a(k)$ is an even function and $b(k)$ is odd in view of the definition \eqref{mtb44a}.

Since $\det[\varphi_\pm,\phi_\pm]=1$, from \eqref{mtb45}, we have
\begin{equation}\label{mtb48}
a(k)=\det[\varphi_-(x;k)\mathrm{e}^{i(k^2-1)x},\phi_+(x;k)\mathrm{e}^{-i(k^2-1)x}]=\det[\varphi_-(0;k),\phi_+(0;k)],
\end{equation}
\begin{equation}\label{mtb48a}
b(k)=\det[\varphi_+(x;k)\mathrm{e}^{-i(k^2-1)x},\varphi_-(x;k)\mathrm{e}^{-i(k^2-1)x}]=\det[\varphi_+(0;k),\varphi_-(0;k)].
\end{equation}
Moreover, the symmetry condition \eqref{mtb47} implies that
\begin{equation}
\label{eq67}
\begin{aligned}
&| a(k)|^{2}+| b(k)|^{2}=1 \quad k\in \mathbb{R},\\
&|a(k)|^{2}-| b(k)|^{2}=1 \quad k\in i\mathbb{R}.
\end{aligned}
\end{equation}

\begin{lemma}\label{lem5}
If $u\in H^{1,1}(\mathbb{R})$, then the function $a(k)$ is continued analytically in $\mathbb{C}^+$ with respect to $z=k^2$. In addition,
\begin{equation}\label{mtb50}
a(k)-1, \quad kb(k), \quad k^{-1}b(k)\in H_z^1(\mathbb{R}),
\end{equation}
 Moreover, if $u\in H^2(\mathbb{R})\cap H^{1,1}(\mathbb{R})$, then
\begin{equation}\label{mtb51}
kb(k), \quad k^{-1}b(k)\in L_z^{2,1}(\mathbb{R}).
\end{equation}

In addition, the maps
\begin{equation}\label{mtb61}
H^{1,1}(\mathbb{R})\ni u\mapsto a(k)-1, \quad kb(k), \quad k^{-1}b(k)\in H_z^1(\mathbb{\mathbb{R}}),
\end{equation}
and
\begin{equation}\label{mtb62}
H^2(\mathbb{R})\cap H^{1,1}(\mathbb{R})\ni u\mapsto kb(k), \quad k^{-1}b(k)\in L_z^{2,1}(\mathbb{R}),
\end{equation}
are Lipshitz continuous.
\end{lemma}
{\bf Proof}~~ Substituting the definition \eqref{mtb44a} into \eqref{mtb43}, and considering the asymptotic behaviors \eqref{mtb10} of the Jost functions $m_-(x;k)$, we have
\begin{equation}\label{mtb57}
\lim\limits_{z\to\infty}a(k)=1.
\end{equation}

For the scattering coefficient $a(k)$, by \eqref{mtb48} and \eqref{mtb44a}, we have
\begin{equation}\label{mtb58}
\begin{aligned}
a(k)-1&=\varphi_-^{(1)}(0;z)\phi_+^{(2)}(0;z)-\varphi_-^{(2)}(0;k)\phi_+^{(1)}(0;k)-1\\
&=[m_-^{(1)}(0;z)-1][n_+^{(2)}(0;z)-1]\\
&+[n_{+}^{(2)}(0;k)-1]+[m_-^{(1)}(0;k)-1]+k^{-1}\varphi_{-}^{(2)}(0;k)k\phi_{+}^{(1)}(0;k).
\end{aligned}
\end{equation}
From the definition \eqref{mtb44a}, we find that
\begin{equation}\label{mtb59}
 2k\phi_+^{(1)}(x;k)+iu(x)
=-n^{(1)}_{+}(x;z)-iu(x)[n^{(2)}_{+}(x;z)-1]\in H_z^1(\mathbb{R}),
\end{equation}
By \eqref{eq44}, we get
\begin{equation}\label{mtb60}
\begin{aligned}
&k^{-1}\varphi_-^{(2)}(x;k)+\frac{i}{2}\int_{-\infty}^x\mathrm{e}^{i(k^2+1)(x-y)}|u(y)|^2k^{-1}\varphi_-^{(2)}(y;k)\mathrm{d}y\\
&=-\int_{-\infty}^x\mathrm{e}^{2i(z+1)(x-y)} \bar{u}(y)[m_-^{(1)}(y,z)-1]\mathrm{d}y
-\int_{-\infty}^x\mathrm{e}^{2i(z+1)(x-y)} \bar{u}(y)\mathrm{d}y,
\end{aligned}
\end{equation}
If we introduce an operator $\mathcal{T}$
\[{\mathcal{T}}f=\frac{i}{2}\int_{-\infty}^x\mathrm{e}^{-2i(k^2-1)(x-y)}|u(y)|^2f(y)\mathrm{d}y,\]
then $(I+{\mathcal{T}})^{-1}$ is exist and bounded for $(\|u\|_{L^2})/2<1$.
Thus, from \eqref{mtb60}, we can get $k^{-1}\varphi_-^{(2)}(x;k)\in H_z^1(\mathbb{R})$. By \eqref{mtb59}, we have $a(k)-1\in H_z^1(\mathbb{R})$.

Next, we consider $kb(k)$ and $k^{-1}b(k)$. By \eqref{mtb48a} and \eqref{mtb44a}, we have
\[\begin{aligned}
2kb(k)=&[m_{-}^{(1)}(0;z)-1]m_+^{(2)}(0;z)+m_+^{(2)}(0;z)\\
&-[m_+^{(1)}(0;z)-1]m_{-}^{(2)}(0;z)-m_{-}^{(2)}(0;z)\in H_z^1(\mathbb{R}),
\end{aligned}\]
and
\[k^{-1}b(k)=m_+^{(1)}(0;z)k^{-1}\varphi_-^{(2)}(0;k)-m_-^{(1)}(0;z)k^{-1}\varphi_+^{(2)}(0;k)\in H_z^1(\mathbb{R}).\]
This finish the proof of \eqref{mtb50}.

Rewritting $2izkb(k)$ as
\[\begin{aligned}
2izkb(k)&=m_-^{(1)}(0;z)[2izm_+^{(2)}(0;z)+\overline{w(0)}]\\
&-m^{(1)}_{+}(0;z)[2izm_-^{(2)}(0;z)+\overline{w(0)}]\\
&+\overline{w(0)}[m^{(1)}_{-}(0;z)-1]\\
&-\overline{w(0)}[m^{(1)}_{+}(0;z)-1]\in L_z^2(\mathbb{R}).
\end{aligned}\]
If $u\in H^2\cap H^{1,1}$, then $kb(k)\in L_z^{2,1}(\mathbb{R})$. $k^{-1}b(k)\in L_z^{2,1}(\mathbb{R})$ can be proved similarly.

To prove the Lipschitz continuity of the maps, we find
\begin{equation}
\label{eq66}
\begin{aligned}
&\|a(k)-\tilde{a}(k)\|_{H^{1}_{z}(\mathbb{R})}\\
&+\| kb(k)-k\tilde{b}(k)\|_{H^{1}_{z}(\mathbb{R})}+\| k^{-1}b(k)-k^{-1}\tilde{b}(k)\|_{H^{1}_{z}(\mathbb{R})}\\
&\leq C(B)\| u-\tilde{u}\|_{H^{1.1}}.
\end{aligned}
\end{equation}
In addition, we have
\begin{equation}
\label{eq66a}
\| kb(k)-k\tilde{b}(k)\|_{L^{2,1}_{z}(\mathbb{R})}+\| k^{-1}b(k)-k^{-1}\tilde{b}(k)\|_{L^{2,1}_{z}(\mathbb{R})}\leq C(B)\| u-\tilde{u}\|_{H^{1,1}(\mathbb{R})\cap H^{2}(\mathbb{R})}.
\end{equation}
Here $C(B)$ is $B$-dependent positive constant. \qquad $\square$

\setcounter{equation}{0}
\section{Riemann-Hilbert problems}\label{sec3}
From \eqref{eq67}, we find that $|a(k)|>1, k\in i\mathbb{R}$. While for $k\in \mathbb{R}$, according to equations \eqref{mtb43} and \eqref{mtb44a}, we find that
\[\begin{aligned}
a(k)-1=\frac{1}{2}\int_{-\infty}^\infty u(x)m_-^{(2)}(x;z)\mathrm{d}x
&:=\tilde{\cal{T}}m_-(x;z)=\tilde{\cal{T}}(I-\tilde{{{K}}}_-)^{-1}e_1.
\end{aligned}
\]
then for every $u\in L^3(\mathbb{R})\cap L^1(\mathbb{R}), u_x\in L^1(\mathbb{R})$, we have
\[|a(k)-1|\leq\frac{1}{2}\|u\|_{L^1}\mathrm{e}^{\|W_1\|_{L^1}},\]
which implies that
\begin{equation}\label{mtc1}
|a(k)|\geq 1-\frac{1}{2}\|u\|_{L^1}\mathrm{e}^{\|W_1\|_{L^1}}>0, \quad k\in\mathbb{R},
\end{equation}
for small enough $u$. Hence, for small enough $u$, $a(k)$ has no zero in $k$( or $z$ in view of that $a(k)$ is even) space.

We define two sectionally analytic functions
\begin{equation}\label{mtc2}
\Psi_+(x;k)=\bigg(\frac{\varphi_-(x;k)}{a(k)},\phi_+(x;k)\bigg), \quad
\Psi_-(x;k)=\bigg(\varphi_+(x;k),\frac{\phi_-(x;k)}{\overline{a(k)}}\bigg),
\end{equation}
where $\Psi_+(x;k)$ analytic in the first and third quadrant of the $k$ plane (that is ${\rm Im} k^2>0$), $\Psi_-(x;k)$ analytic in the second and fourth quadrant of the $k$ plane (where ${\rm Im} k^2<0$), and
\begin{equation}\label{mtc7}
\Psi_\pm(x;k)\to I, \quad k\to\infty,
\end{equation}
 For every $x\in\mathbb{R}$, there exists a jump condition
\begin{equation}\label{mtc3}
 \Psi_+(x;k)-\Psi_-(x;k)=\Psi_-(x;k)S(x;k), \quad k\in\mathbb{R}\cup i\mathbb{R},
\end{equation}
where
\begin{equation}
\label{mtc4}
S(x;k)=\left(\begin{matrix}
\mid r(k)\mid^{2} &  \overline{r(k)}\mathrm{e}^{-2 i(k^{2}+1)x}\\
r(k)\mathrm{e}^{ 2i(k^{2}+1)x}  &  0
\end{matrix}\right) ,\quad k\in \mathbb{R},
\end{equation}
\begin{equation}
\label{mtc5}
S(x;k)=\left(\begin{matrix}
-\mid r(k)\mid^{2} & - \overline{r(k)}\mathrm{e}^{ -2i(k^{2}+1)x}\\
r(k)\mathrm{e}^{2i(k^{2}+1)x}  &  0
\end{matrix}\right) ,\quad k\in  i\mathbb{R},
\end{equation}
with the reflection coefficient
\begin{equation}\label{mtc6}
 r(k)=\frac{b(k)}{a(k)}.
\end{equation}
Thus, \eqref{mtc2},\eqref{mtc3} and \eqref{mtc7} construct a Riemann-Hilbbert problem. When $k\in \mathbb{R}$, matrix $S$ is Hermitian, then the Riemann-Hilbert problem \eqref{mtc3},\eqref{mtc5} and \eqref{mtc7} exist a unique solution. However, when $k\in i\mathbb{R}$ the matrix $S$ is not a Hermitian. But from \eqref{eq67}, we find that
\begin{equation}\label{mtc8}
 1-|r(k)|^2=\frac{1}{|a(k)|^2}\geq C_0^2>0, \quad k\in i\mathbb{R},
\end{equation}
where $C_0=\sup\limits_{k\in\mathbb{R}}|a(k)|$. The condition \eqref{mtc8} implies that the Riemann-Hilbert problem \eqref{mtc3},\eqref{mtc4} and \eqref{mtc7} also exist a unique solution.

Now we define two new reflection coefficients
\begin{equation}\label{mtc9}
 r_+(z)=\frac{r(k)}{2k}, \quad r_-(z)=2kr(k),
\end{equation}
which satisfy
\begin{equation}\label{mtc10}
r_-(z)=4zr_+(z),
\end{equation}
and
\begin{equation}\label{mtc11}
\begin{aligned}
\overline{r_+(z)}r_-(z)=|r(k)|^2, \quad z\in\mathbb{R}^+,k\in\mathbb{R},\\
\overline{r_+(z)}r_-(z)=-|r(k)|^2, \quad z\in\mathbb{R}^-,k\in i\mathbb{R}.
\end{aligned}
\end{equation}

\begin{lemma}\label{lem8}
Assume $a(k)$ satisfies the condition \eqref{mtc1}. If $u\in H^{1,1}(\mathbb{R})$, then $r_\pm(z)\in H^1(\mathbb{R})$,
and if  $u\in H^2(\mathbb{R})\cap H^{1,1}(\mathbb{R})$, then $r_\pm(z)\in L^{2,1}(\mathbb{R})$. Moreover, the map
\begin{equation}\label{mtc12}
 H^2(\mathbb{R})\cap H^{1,1}(\mathbb{R})\ni u\mapsto(r_-,r_+)\in H^{1,1}(\mathbb{R})\cap L^{2,1}(\mathbb{R}),
\end{equation}
is Lipschitz continuous.
\end{lemma}
{\bf Proof}~~ According to Lemma \ref{lem5}, it can be inferred that if $u\in H^{1,1}(\mathbb{R})$, then $k^{-1}b(k)\in H^1(\mathbb{R})$, so $r_\pm(z)\in H^1(\mathbb{R})$. We note that
\begin{equation}\label{mtc13}
\begin{aligned}
r_{-}-\tilde{r}_{-}=\frac{2k(b-\tilde{b})}{a}+\frac{2k\tilde{b}}{a\tilde{a}}[(\tilde{a}-1)-(a-1)],\\
r_{+}-\tilde{r}_{+}=\frac{(b-\tilde{b})}{2ka}+\frac{\tilde{b}}{2ka\tilde{a}}[(\tilde{a}-1)-(a-1)].
\end{aligned}
\end{equation}
Accordingly, the Lipschitz continuity can be proved. \qquad $\square$

Next, we define two new sectionally analytic functions
\begin{equation}\label{mtc14}
\begin{aligned}
M_+(x;z)=\left(\frac{m_-(x;z)}{a(z)},p_+(x;z)\right), \\ M_-(x;z)=\left(m_+(x;z),\frac{p_-(x;z)}{\overline{a(z)}}\right),
\end{aligned}
\end{equation}
where
\begin{equation}\label{mtc15}
 p_\pm(x;z)=\frac{1}{2k}T_1(x;k)T_2^{-1}(x;k)n_\pm(x;z).
\end{equation}
Note that $a(-k)=a(k)$, we write it as $a(z), z=k^2$. Then the jump condition \eqref{mtc3} can be rewritten as
\begin{equation}\label{mtc16}
M_+(x;z)-M_-(x;z)=M_-(x;z)R(x;z), \quad z\in \mathbb{R},
\end{equation}
where
\begin{equation}\label{mtc16a}
R(x;z)=\left(
\begin{matrix}
\overline{r_+(z)}r_-(z)&\overline{r_+(z)}\mathrm{e}^{-2i(z+1)x}\\
r_-(z)\mathrm{e}^{2i(z+1)x}&0
\end{matrix}\right).
\end{equation}

Since equation \eqref{mtc15} satisfy
$\lim\limits_{z\to\infty}p_\pm(x,;z)=e_2$,
then
\begin{equation}\label{mtc18}
 M_\pm(x;z)\to I, \quad z\to\infty.
\end{equation}

For the Riemann-Hilbert problem \eqref{mtc14},\eqref{mtc16} and \eqref{mtc18}, the matrix $R(x;z)$ is not Hermitian. So it is difficult to construct its unique solution. Here, for every $k\in \mathbb{C}\setminus\{0\}$, we introduce two matrices
\begin{equation}\label{mtc19}
\tau_1(k)={\rm diag}(1,2k), \quad \tau_2(k)=(2k)^{-1}\tau_1(k),
\end{equation}
then we find that
\[\tau_j^{-1}(k)R(x;z)\tau_j(k)=S(x;k), \quad j=1,2.\]
Furthermore, the jump equation \eqref{mtc16} can be rewritten as
\begin{equation}\label{mtc20}
 \begin{aligned}
\Xi_{+1,2}(x;k)-\Xi_{-1,2}(x;k)&=\Xi_{-1,2}(x;k)S(x;k)+F_{1,2}(x;k), \quad k\in\mathbb{R}\cap i\mathbb{R},\\
\Xi_{\pm1,2}(x;k)&\to0, \quad k\to\infty,
 \end{aligned}
\end{equation}
where
\begin{equation}\label{mtc21}
 \Xi_{\pm1,2}(x;k)=(M_\pm(x;z)-I)\tau_{1,2}(k), \quad F_{1,2}(x;k)=\tau_{1,2}(k)S(x;k).
\end{equation}
and  $\Xi_{-1,2}(x;\cdot)$ is analytic in the first and third quadrant of the $k$ plane, and $\Xi_{+1,2}(x;\cdot)$ is analytic in the second and fourth quadrant of the $k$ plane.

\setcounter{equation}{0}
\section{The inverse scattering transform and Lipschitz continuity}\label{sec4}
\subsection{Solvability of the Riemann-Hilbert problems}
Under the constraint \eqref{mtc12}, we will discuss the solvability of the Riemann-Hilbert problem \eqref{mtc20}.

\begin{proposition}\label{pro2}
If $r_\pm(z)\in H^1(\mathbb{R})\cap L^{2,1}(\mathbb{R})$, then $r(k)\in L_z^{2,1}(\mathbb{R})\cap L_z^\infty(\mathbb{R})$.
In addition, if $r_-(z)\in H^1(\mathbb{R})\cap L^{2,1}(\mathbb{R})$, then $\|kr_-(z)\|_{L_z^\infty}\leq\|r_-\|_{H^1\cap L^{2,1}(\mathbb{R})}$
\end{proposition}
{\bf Proof}.~~From \eqref{mtc11}, we know that $|r(k)|^2=-{\rm sgn }z\overline{r_+(z)}r_-(z)$. So, by virtue of the Schwartz inequality, we find that if $r_\pm(z)\in L^{2,1}(\mathbb{R})$, then $r(k)\in L^{2,1}(\mathbb{R})$.

Since $\|f\|_{L^\infty}\leq\frac{1}{\sqrt{2}}\|f\|_{H^1}$, then $r_\pm\in L_z^\infty$ for $r_\pm\in H^1(\mathbb{R})$ and furthermore $r(k)\in L_z^\infty(\mathbb{R})$.

To prove the second statement, we find that
\[\begin{aligned}
|zr_-(z)|^2&\leq\int_0^z|r_-^2(z)+2zr_-(z)\partial_zr_{-}(z)|\mathrm{d}z\\
&\leq\|r_-\|_{L^2}^2+2\|r_-\|_{L^{2,1}}\|\partial_zr_-\|_{L_z^{2}}\leq\|r_-\|^2_{H^1\cap L^{2,1}(\mathbb{R})}.
\end{aligned}\]
The Proposition is proved.  \qquad $\square$

Now, we introduce the Cauchy operator and the project operator. For every $f\in L^p(\mathbb{R}), 1\leq p<\infty$, the Cauchy operator ${\mathcal{C}}_z$ is
\begin{equation}
\label{mtd1}
{\mathcal{C}}_z(f)=\frac{1}{2\pi i}\int_{\mathbb{R}}\frac{f(s)}{s-z}\mathrm{d}s,\quad z\in \mathbb{C}\setminus \mathbb{R}.
\end{equation}
The projection operators ${\mathcal{P}}^{\pm}$ are given explicitly
\begin{equation}\label{mtd2}
{\mathcal{P}}^{\pm}(f)(z)=\frac{1}{2\pi i}\lim\limits_{\varepsilon\to 0}\int_{\mathbb{R}}\frac{f(s)}{s-(z\pm\varepsilon i)}\mathrm{d}s,\quad z\in \mathbb{R}.
\end{equation}

\begin{lemma}\cite{imrn2018-5663}\label{lem9}
For every $r(k)\in L_z^2(\mathbb{R})\cap L_z^\infty(\mathbb{R})$ satisfying the constraint \eqref{mtc8}, and for every $F(k)\in L_z^2(\mathbb{R})$, there exist a unique solution $G(k)\in L_z^2(\mathbb{R})$ to the linear equation
\begin{equation}\label{mtd3}
 (I-{\mathcal{P}}^{-}_{s})G(k)=F(k),\quad k\in \mathbb{R}\cup i\mathbb{R},
\end{equation}
where ${\mathcal{P}}^{-}_{s}G={\mathcal{P}}^{-}(GS)$. In addition, the operator $(I-{\mathcal{P}}^{-}_{s})$ is inverse and for every $f(z)\in L^2(R)$ ,
\begin{equation} \label{mtd4}
 \|(I-{\mathcal{P}}^{-}_{s})^{-1}f\|_{{L}^{2}_{z}(\mathbb{R})}\leq C\| f\|_{{L}^{2}_{z}(\mathbb{R})},
\end{equation}
where $C$ is a constant depending on $\| r(k)\|_{L^{\infty}_{z}(\mathbb{R})}$.
\end{lemma}

\begin{corollary}\label{cor1}
Let $r_\pm(z)\in H^1(\mathbb{R})\cap L^{2,1}(\mathbb{R})$ satisfying the constraint \eqref{mtc10} and \eqref{mtc8}. For every $x\in\mathbb{R}$,
the Riemann-Hilbert problem given in \eqref{mtc20} exist a unique solution
\[\Xi_{\pm1,2}(x;k)=M_{\pm}(x;z)\tau_{1,2}(k)-\tau_{1,2}(k).\]
\end{corollary}
{\bf Proof}.~~ From \eqref{mtd3}, we find the unique solution $\Xi_{-1,2}(x;k)$ to the equation
\begin{equation}\label{mtd5}
\Xi_{-1,2}(x;k)={\mathcal{P}}^{-}(\Xi_{-1,2}(x;k)S(x;k)+F_{1,2}(x;k))(z),\quad z\in\mathbb{R}.
\end{equation}
We define $\Xi_{+1,2}(x;k)$ by
\begin{equation}\label{mtd6}
\Xi_{+1,2}(x;k)={\mathcal{P}}^{+}(\Xi_{-1,2}(x;k)S(x;k)+F_{1,2}(x;k))(z),\quad z\in\mathbb{R}.
\end{equation}
Then the unique solution of the Riemann-Hilbert problem can be obtained by analytic continuation to the Cauchy operator
\begin{equation}
\label{mtd7}
\Xi_{\pm1,2}(x;k)={\mathcal{C}}_z\big(\Xi_{-1,2}(x;k)S(x;k)+F_{1,2}(x;k)\big),\quad z\in \mathbb{C}^{\pm}.
\end{equation}
The corollary is proved. \qquad $\square$

Next, we will provide the estimates of solutions to the Riemann-Hilbert problem given in \eqref{mtc16}.
We define  $M_\pm$ as
\begin{equation}\label{mtd8}
M_\pm(x;z)=\big(\xi_\pm(x;z), \eta_\pm(x;z)\big).
\end{equation}
Then from \eqref{mtc21}
\begin{equation}\label{mtd9}
\Xi_{\pm1}(x;k)=M_{\pm}(x;z)\tau_{1}(k)-\tau_{1}(k)=\big(\xi_{\pm}-e_{1},2k(\eta_{\pm}-e_{2})\big),
\end{equation}
and
\begin{equation}\label{mtd10}
\Xi_{\pm2}(x;k)=M_{\pm}(x;z)\tau_{2}(k)-\tau_{2}(k)=\big((2k)^{-1}(\xi_{\pm}-e_{1}),\eta_{\pm}-e_{2}\big).
\end{equation}
Since
\[F_{1,2}(x;k)=\tau_{1,2}(k)S(x;k)=R(x;z)\tau_{1,2}(k),\]
then
\begin{equation}\label{mtd11}
\Xi_{\pm1,2}(x;k)S(x;k)+F_{1,2}(x;k)=M_{\pm}(x;z)R(x;z)\tau_{1,2}(k).
\end{equation}
Substituting it into \eqref{mtd5},\eqref{mtd6} and taking the first column for the solution $\Xi_{\pm1}(x;k)$, we have
\begin{equation}\label{mtd12}
\xi_{\pm}-e_{1}={\mathcal{P}}^{\pm}((M_{-}R)^{[1]}(x;\cdot))(z),
\end{equation}
Similarly
\begin{equation}\label{mtd13}
\eta_{\pm}-e_{2}={\mathcal{P}}^{\pm}((M_{-}R)^{[2]}(x;\cdot))(z).
\end{equation}
\begin{equation}\label{mtd12a}
\begin{aligned}
(2k)^{-1}(\xi_{\pm}-e_{1})&={\mathcal{P}}^{\pm}((2k)^{-1}(M_{-}R)^{[1]}(x;\cdot))(z),\\
2k(\eta_{\pm}-e_{2})&={\mathcal{P}}^{\pm}((2k)(M_{-}R)^{[2]}(x;\cdot))(z).
\end{aligned}
\end{equation}

By combining equations \eqref{mtd12} and \eqref{mtd13}, we can obtain
\begin{equation}\label{mtd14}
M_{\pm}(x;z)=I+{\mathcal{P}}^{\pm}(M_{-}(x;\cdot)R(x;\cdot))(z),\quad z\in \mathbb{R},
\end{equation}
and the solution of the Riemann-Hilbert problem \eqref{mtc16}
\begin{equation}\label{mtd15}
M_{\pm}(x;z)=I+{\mathcal{C}}(M_{-}(x;\cdot)R(x;\cdot))(z),\quad z\in \mathbb{C}^\pm.
\end{equation}

\begin{lemma}\label{lem10}
Let $r_\pm(z)\in H^1(\mathbb{R})\cap L^{2,1}(\mathbb{R})$ satisfying the constraint \eqref{mtc10} and \eqref{mtc8}.
The unique solution to the
integral equations \eqref{mtd14} satisfies the estimate
\begin{equation}\label{mtd16}
\|M_{\pm}(x;\cdot)-I\|_{{L}^{2}(\mathbb{R})}\leq C(\| r_{+}\|_{{L}^{2}(\mathbb{R})}+\| r_{-}\|_{{L}^{2}(\mathbb{R})}).
\end{equation}
\end{lemma}
Where $C$ is a positive constant that only depending on $\|r_\pm\|_{L^\infty}$.

{\bf Proof}.~~ From Proposition \ref{pro2}, we know that if $r_{\pm}(z)\in {H}^{1}(\mathbb{R})\cap {L}^{2.1}(\mathbb{R})$, then $r(k)\in {L}^{2}_{z}(\mathbb{R})\cap {L}^{\infty}_{z}(\mathbb{R})$ and
\begin{equation}\label{mtd17}
\| R(x;z)\tau_{1,2}(k)\| _{{L}^{2}_{z}(\mathbb{R})}\leq C(\| r_{+}\|_{{L}^{2}(\mathbb{R})}+\|r_{-}\|_{{L}^{2}(\mathbb{R})}).
\end{equation}
From \eqref{mtd12} and \eqref{mtd9}, we know that, for every $x\in\mathbb{R}$
\[\begin{aligned}\|(M_{\pm}(x;z)-I)^{[1]}\|_{{L_z}^{2}(\mathbb{R})}&=\|{\mathcal{P}}^{\pm}(M_{-}(x;k)R(x;k))^{[1]}(z)\|_{L_z^2}\\
&=\|\Xi_{\pm1}^{[1]}(x;k)\|_{L_z^2}\\
&\leq C_1\|(R(x;k)\tau_1(k))^{[1]}\|_{L_z^2}\\
&\leq C(\| r_{+}\|_{{L}^{2}(\mathbb{R})}+\|r_{-}\|_{{L}^{2}(\mathbb{R})}),\end{aligned}\]
where $C,C_1$ is a positive constant only depend on $\|r\|_{L^\infty}$.
The proof for $(M_{\pm}(x;z)-I)^{[2]}$ is analogous.

Using \eqref{mtc14} and \eqref{mtc16},
from \eqref{mtd12}, \eqref{mtd13} and \eqref{mtd8}, we have
\begin{equation}\label{mtd18}
\begin{aligned}
\xi_-(x;z)-e_1={\cal{P}}^-\left(r_-(z)\mathrm{e}^{2i(z+1)x}\eta_+(x;z)\right),\\
\eta_+(x;z)-e_2={\cal{P}}^+\left(\overline{r_+(z)}\mathrm{e}^{-2i(z+1)x}\xi_-(x;z)\right).
\end{aligned}
\end{equation}
Defining a new matrix
\begin{equation}\label{mtd19}
 M(x;z)=(\xi_-(x;z)-e_1,\eta_+(x;z)-e_2),
\end{equation}
then \eqref{mtd18} can be rewritten as
\begin{equation}
\label{b64}
M(x;z)-{\mathcal{P}}^{+}(M(x;z)R_{+}(x;z))-{\mathcal{P}}^{-}(M(x;z)R_{-}(x;z))=F(x;z),
\end{equation}

Since ${\mathcal{P}}^{+}-{\mathcal{P}}^{-}=I$ and $R_{+}(x;z)+R_{-}(x;z)=\big(I-R_{+}(x;z)\big)R(x;z),$
then equation \eqref{b64} can be rewritten as
\begin{equation}\label{mtd23}
G(x;z)-{\mathcal{P}}^{-}(G(x;z)R(x;z))=F(x;z),
\end{equation}
where
\begin{equation}\label{mtd24}
\begin{aligned}
G(x;z)&=M(x;z)\big(I-R_+(x;z)\big)\\
&=\left(\begin{matrix}
\xi_-^{(1)}-1&\eta_+^{(1)}-\overline{r_+}\mathrm{e}^{-2i(z+1)x}(\xi_-^{(1)}-1)\\
\xi_-^{(2)}&\eta_+^{(2)}-1-\overline{r_+}\mathrm{e}^{-2i(z+1)x}\xi_-^{(2)}
\end{matrix}\right).
\end{aligned}
\end{equation}
\begin{equation}
\label{b65}
R_{+}(x;z)=\left(\begin{matrix}
0 &  \overline{r_{+}(z)}\mathrm{e}^{-2i(z+1)x}\\
0   &  0
\end{matrix}\right),\quad R_{-}(x;z)=\left(\begin{matrix}
0 &  0\\
r_{-}(z)\mathrm{e}^{2i(z+1)x}   &  0
\end{matrix}\right),
\end{equation}
\begin{equation}
\label{b66}
F(x;z)=\big(e_{2}{\mathcal{P}}^{-}(r_{-}(z)\mathrm{e}^{2i(z+1)x}),
e_{1}{\mathcal{P}}^{+}(\overline{r_{+}(z)}\mathrm{e}^{-2i(z+1)x})\big).
\end{equation}
Note that $R(x;z)\tau_{1,2}(k)=\tau_{1,2}(k)S(x;k)$, equation \eqref{mtd23} can be rewritten as
\begin{equation}\label{mtd35}
 G_{1,2}(x;k)-{\cal{P}}_s^-(G_{1,2}(x;k))=F(x;z)\tau_{1,2}(k),
\end{equation}
where $G_{1,2}(x;k)=G(x;z)\tau_{1,2}(k)$ and the operator ${\cal{P}}_s^-$ is defined in Lemma \ref{lem9}.
Since the second row of $F(x;z)$ and $F(x;z)\tau_1(k)$ share the same vector $({\mathcal{P}}^{-} (r_{-}(z)\mathrm{e}^{2i(z+1)x}),0)$, then for the case of $G_1(x;k)$ in \eqref{mtd35}, considering its second row and using the bound \eqref{mtd4}, we have
\begin{equation}
\label{b54}
\| \xi^{(2)}_{-}(x;z)\|_{{L}^{2}_{z}(\mathbb{R})}\leq C_1\|{\mathcal{P}}^{-}(r_{-}(z)\mathrm{e}^{2i(z+1) x})\|_{{L}^{2}_{z}(\mathbb{R})}.
\end{equation}

Similarly, the first row of $F(x;z)$ and $F(x;z)\tau_2(k)$ share the same vector $(0,{\mathcal{P}}^{+}(\overline{r_{+}(z)}\mathrm{e}^{-2i(z+1)x}))$,
then for the case of $G_2(x;k)$ in \eqref{mtd35}, considering its first row and using the bound \eqref{mtd4}, we have
\begin{equation}
\label{b71}
\| (2k)^{-1}(\xi^{(1)}_{-}-1)\|_{{L}^{2}_{z}(\mathbb{R})}\leq C\|{\mathcal{P}}^{+}(\overline{r_{+}(z)}\mathrm{e}^{-2i(z+1)x})\|_{{L}^{2}_{z}(\mathbb{R})},
\end{equation}
\begin{equation}
\label{b72}
\| \eta^{(1)}_{+}(x;z)-\overline{r}_{+}(z)\mathrm{e}^{-2i(z+1)x}(\xi^{(1)}_{-}-1)\|_{{L}^{2}_{z}(\mathbb{R})}\leq C\|{\mathcal{P}}^{+}(\overline{r_{+}(z)}\mathrm{e}^{-2i(z+1)x})\|_{{L}^{2}_{z}(\mathbb{R})},
\end{equation}
which imply that
\begin{equation}
\label{b73}
\begin{aligned}
\| \eta^{(1)}_{+}(x;z)\|_{L^{2}_{z}}
&\leq C_2\|{\mathcal{P}}^{+}(\overline{r_{+}(z)}\mathrm{e}^{-2i(z+1)x}\|_{{L}^{2}_{z}(\mathbb{R})},
\end{aligned}
\end{equation}
in terms of $|2k\overline{r_{+}(z)}|=| r(k)|$, where $C_2$ are a positive constant depending on $\| r_{\pm}\|_{L^{\infty}(\mathbb{R})}$.

We note that
\begin{equation}\label{b55}
\begin{aligned}
\|\partial_x \xi^{(2)}_{-}(x;z)\|_{{L}^{2}_{z}(\mathbb{R})}\leq& C_3\big(\|\partial_x{\mathcal{P}}^{-}(r_{-}(z)\mathrm{e}^{-2i(z+1)x})\|_{{L}^{2}_{z}(\mathbb{R})}\\
&\quad+\|{\mathcal{P}}^{-}(r_{-}(z)\mathrm{e}^{-2i(z+1)x})\|_{{L}^{2}_{z}(\mathbb{R})}\big),
\end{aligned}
\end{equation}
and
\begin{equation}\label{mtd38}
\begin{aligned}
\|\partial_x \eta^{(1)}_{+}(x;z)\|_{L^{2}_{z}}\leq& C_4\big(\|\partial_x{\mathcal{P}}^{+}(\overline{r_{+}(z)}\mathrm{e}^{-2i(z+1)x}\|_{{L}^{2}_{z}(\mathbb{R})}\\
&\quad+\|{\mathcal{P}}^{+}(\overline{r_{+}(z)}\mathrm{e}^{-2i(z+1)x}\|_{{L}^{2}_{z}(\mathbb{R})},
\end{aligned}
\end{equation}
which can be proved by differentiating equation \eqref{b64} with respect to $x$.

\subsection{Reconstruction for the potential \texorpdfstring{$u$}{}}\label{sec5}
Now let's discuss the potential reconstruction.
From \eqref{mtb12}, we have
\begin{equation}\label{mte22}
i\overline{u_x(x)}+2\overline{u}-\frac{1}{2}|u|^{2}\overline{u}=-2i\lim\limits_{z\to\infty}zm_\pm^{(2)}(x;z).
\end{equation}
In addition, from \eqref{mtc15}, we have
\begin{equation}\label{mte23}
iu(x)=-4\lim\limits_{z\to\infty}zp_\pm^{(1)}(x;z).
\end{equation}

Firstly, we consider the case of $x\in\mathbb{R}^+$.
Using the expression of the solution \eqref{mtd15}, we rewrite the reconstraction formulas \eqref{mte22} and \eqref{mte23}
\begin{equation}\label{mte27}
\begin{aligned}
i\overline{u_x(x)}+2\overline{u}-\frac{1}{2}|u|^{2}\overline{u}=\frac{1}{\pi}\int_{\mathbb{R}}r_-(z)\mathrm{e}^{2i(z+1)x}\eta_+^{(2)}(x;z)\mathrm{d}z,
\end{aligned}
\end{equation}
and
\begin{equation}\label{mte28}
iu(x)=\frac{2}{i\pi}\int_{\mathbb{R}}\overline{r_+(z)}\mathrm{e}^{-2i(z+1)x}\xi_-^{(1)}(x;z)\mathrm{d}z.
\end{equation}
From the integral equations \eqref{mte27} and \eqref{mte28}, we can obtain the following estimate.
\begin{lemma}\label{lem12}
If $r_{\pm}(z)\in {H}^{1}(\mathbb{R})\cap {L}^{2,1}(\mathbb{R})$ and admit the condition \eqref{mtc8}, then $u\in{H}^{2}(\mathbb{R}^+)\cap {H}^{1,1}(\mathbb{R}^+)$. In addition, we have
\begin{equation}
\label{b74}
\| u\|_{{H}^{2}(\mathbb{R}^{+})\cap {H}^{1,1}(\mathbb{R}^{+})}\leq C(\| r_{+}\|_{{H}^{1}(\mathbb{R})\cap {L}^{2,1}(\mathbb{R})}+\| r_{-}\|_{{H}^{1}(\mathbb{R})\cap {L}^{2,1}(\mathbb{R})}),
\end{equation}
where $C$ is a positive constant depending on $\|r_\pm(z)\|_{L^\infty(\mathbb{R})}$ and $\| r_{\pm}\|_{{H}^{1}(\mathbb{R})\cap {L}^{2,1}(\mathbb{R})}$.
\end{lemma}
{\bf Proof}~~ We rewrite the construction formula \eqref{mte28} as
\begin{equation}\label{mte29}
\begin{aligned}
iu(x)=A_1(x)+A_2(x), \quad
i\overline{u_x(x)}+2\overline{u}-\frac{1}{2}|u|^{2}\overline{u}=B_1(x)+B_2(x),
\end{aligned}
\end{equation}
where
\begin{equation}\label{mte29a}
\begin{aligned}
A_1(x)&=\frac{2}{\pi i}\int_{\mathbb{R}}\overline{r_{+}(z)}\mathrm{e}^{-2i(z+1)x}dz\\
A_2(x)&=\frac{2}{\pi i}\int_{\mathbb{R}}\overline{r_{+}(z)}\mathrm{e}^{-2i(z+1)x}\big(\xi_{-}^{(1)}(x;z)-1\big)\mathrm{d}z,
\end{aligned}
\end{equation}
and
\begin{equation}\label{mte29b}
\begin{aligned}
B_1(x)&=\frac{1}{\pi}\int_{\mathbb{R}}r_-(z)\mathrm{e}^{2i(z+1)x}\mathrm{d}z,\\
B_2(x)&=\frac{1}{\pi}\int_{\mathbb{R}}r_-(z)\mathrm{e}^{2i(z+1)x}\big(\eta_+^{(2)}(x;z)-1\big)\mathrm{d}z.
\end{aligned}
\end{equation}
Using the theories of the Fourier transform, we find that
\begin{equation}\label{mtea3}
\begin{aligned}
\|A_1(x)\|_{L^2(\mathbb{R}^+)}=\frac{2}{\sqrt{\pi}}\|r_+(z)\|_{L^2(\mathbb{R})}, \\
\|B_1(x)\|_{L^2(\mathbb{R}^+)}=\frac{1}{\sqrt{\pi}}\|r_-(z)\|_{L^2(\mathbb{R})}.
\end{aligned}
\end{equation}
In addition, the differential properties of the Fourier transform imply that, for every $x\in\mathbb{R}^+$
\begin{equation}\label{mtea4}
\begin{aligned}
 \|\langle x\rangle A_1(x)\|_{L^2(\mathbb{R}^+)}=\frac{1}{\sqrt{\pi}}\|\partial_zr_+(z)\|_{L^2(\mathbb{R})}, \\
\|\langle x\rangle B_1(x)\|_{L^2(\mathbb{R}^+)}=\frac{1}{2\sqrt{\pi}}\|\partial_zr_-(z)\|_{L^2(\mathbb{R})},
\end{aligned}
\end{equation}
and
\begin{equation}\label{mtea5}
\|\partial_x B_1(x)\|_{L^2(\mathbb{R}^+)}\leq\frac{2}{\sqrt{\pi}}\big(\|r_-(z)\|_{L^{2,1}(\mathbb{R})}+\|r_-(z)\|_{L^2(\mathbb{R})}\big).
\end{equation}

For $A_2(x)$ and $B_2(x)$, using \eqref{mtd18} we can find that
\begin{equation}\label{mteaa5}
\begin{aligned}
A_2(x)&=\frac{2}{i\pi}\int_{\mathbb{R}}\overline{r_{+}(z)}\mathrm{e}^{-2i(z+1)x}{\cal{P}}^-\big(r_-(z)\mathrm{e}^{2i(z+1)x}\eta_+^{(1)}(x;z)\big)\mathrm{d}z\\
&=-\frac{2}{i\pi}\int_{\mathbb{R}}{\cal{P}}^+\big(\overline{r_{+}(z)}\mathrm{e}^{-2i(z+1)x}\big)r_-(z)\mathrm{e}^{2i(z+1)x}\eta_+^{(1)}(x;z)\mathrm{d}z,
\end{aligned}
\end{equation}
and
\begin{equation}\label{mtea6}
B_2(x)=-\frac{1}{\pi}\int_{\mathbb{R}}{\cal{P}}^-\big(\overline{r_{-}(z)}\mathrm{e}^{2i(z+1)x}\big)\overline{r_+(z)}\mathrm{e}^{-2i(z+1)x}\xi_-^{(2)}(x;z)\mathrm{d}z.
\end{equation}
In addition,
\begin{equation}\label{mtea10}
\begin{aligned}
\partial_xB_2(x)=&-\frac{1}{\pi}\int_{\mathbb{R}}[\partial_x{\cal{P}}^-\big(\overline{r_{-}(z)}\mathrm{e}^{2i(z+1)x}\big)]\overline{r_+(z)}\mathrm{e}^{-2i(z+1)x}\xi_-^{(2)}(x;z)\mathrm{d}z\\
&+\frac{2i}{\pi}\int_{\mathbb{R}}{\cal{P}}^-\big(\overline{r_{-}(z)}\mathrm{e}^{2i(z+1)x}\big)(z+1)\overline{r_+(z)}\mathrm{e}^{-2i(z+1)x}\xi_-^{(2)}(x;z)\mathrm{d}z\\
&-\frac{1}{\pi}\int_{\mathbb{R}}{\cal{P}}^-\big(\overline{r_{-}(z)}\mathrm{e}^{2i(z+1)x}\big)\overline{r_+(z)}\mathrm{e}^{-2i(z+1)x}\partial_x\xi_-^{(2)}(x;z)\mathrm{d}z.
\end{aligned}
\end{equation}

By virtue of the inequalities \eqref{b54} and \eqref{b73}, we find that
\begin{equation}\label{nlsgic1}
\begin{aligned}
|A_2(x)|\lesssim\|r_-(z)\|_{L^\infty}\|{\cal{P}}^+\big(\overline{r_{+}(z)}\mathrm{e}^{-2i(z+1)x}\big)\|^2_{L_z^2(\mathbb{R})},\\
|B_2(x)|\lesssim\|r_+(z)\|_{L^\infty}\|{\cal{P}}^-\big({r_{-}(z)}\mathrm{e}^{2i(z+1)x}\big)\|^2_{L_z^2(\mathbb{R})},
\end{aligned}
\end{equation}
and
\begin{equation}\label{nlsgic2}
\begin{aligned}
|\partial_xB_2(x)|\lesssim&\|r_+(z)\|_{L^\infty}\|{\cal{P}}^-\big({r_{-}(z)}\mathrm{e}^{2i(z+1)x}\big)\|_{L_z^2(\mathbb{R})}\|\partial_x{\cal{P}}^-\big({r_{-}(z)}\mathrm{e}^{2i(z+1)x}\big)\|_{L_z^2(\mathbb{R})}\\
&+\big(\|\|r_+(z)\|_{L^\infty}+\|r_-(z)\|_{L^\infty}\big)\|{\cal{P}}^-\big({r_{-}(z)}\mathrm{e}^{2i(z+1)x}\big)\|^2_{L_z^2(\mathbb{R})}.
\end{aligned}
\end{equation}
Here we have used the relation $r_-(z)=4zr_+(z)$.

It is noted that, for $x>0$ \cite{cpam51-697}
\begin{equation}\label{nlsgic3}
\begin{aligned}
\|{\cal{P}}^+\big(\overline{r_{+}(z)}\mathrm{e}^{-2i(z+1)x}\big)\|_{L_z^2(\mathbb{R})}\leq\frac{1}{\sqrt{1+x^2}}\|r_+(z)\|_{H^1(\mathbb{R})},\\
\|{\cal{P}}^-\big({r_{-}(z)}\mathrm{e}^{2i(z+1)x}\big)\|_{L_z^2(\mathbb{R})}\leq\frac{1}{\sqrt{1+x^2}}\|r_+(z)\|_{H^1(\mathbb{R})},
\end{aligned}
\end{equation}
and
\begin{equation}\label{nlsgic4}
\|\partial_x{\cal{P}}^-\big({r_{-}(z)}\mathrm{e}^{2i(z+1)x}\big)\|_{L_z^2(\mathbb{R})}\leq\frac{1}{\sqrt{1+x^2}}\big(\|r_+(z)\|_{H^1(\mathbb{R})}+\|zr_+(z)\|_{H^1(\mathbb{R})}\big).
\end{equation}

Substituting \eqref{nlsgic3}, \eqref{nlsgic4} into \eqref{nlsgic1}, \eqref{nlsgic2}, we can get that
\begin{equation}\label{mtea33}
\begin{aligned}
\|A_2(x)\|^2_{L^2(\mathbb{R}^+)}\lesssim \|r_-(z)\|_{L^\infty}\|r_+(z)\|^2_{H^{1}(\mathbb{R})}),\\
\|B_2(x)\|^2_{L^2(\mathbb{R}^+)}\lesssim \|r_+(z)\|_{L^\infty}\|r_-(z)\|^2_{H^{1}(\mathbb{R})}),
\end{aligned}
\end{equation}
and
\begin{equation}\label{mtea34}
\begin{aligned}
\|\partial_xB_2(x)\|^2_{L^2(\mathbb{R}^+)}\lesssim&(\|r_+(z)\|^2_{L^\infty(\mathbb{R})}+\|r_-(z)\|_{L^\infty(\mathbb{R})})\|r_-(z)\|^2_{H^1\cap L^{2,1}(\mathbb{R})}.
\end{aligned}
\end{equation}

In addition, according to the weighted H\"older inequality and sharp Gagliardo-Nirenberg inequality, we have
\begin{equation}\label{mtea19}
\begin{aligned}
\|\langle x\rangle A_2(x)\|_{L^2(\mathbb{R}^+)}&\lesssim\|r_-(z)\|_{L^\infty(\mathbb{R})}\|\langle x\rangle{\cal{P}}^+\big(\overline{r_{+}(z)}\mathrm{e}^{-2i(z+1)x}\big)\|_{L^\infty(\mathbb{R}^+,L_z^2(\mathbb{R}))}\\
&\qquad\times\|{\cal{P}}^+\big(\overline{r_{+}(z)}\mathrm{e}^{-2i(z+1)x}\big)\|_{L^2(\mathbb{R}^+,L_z^2(\mathbb{R}))}\\
&\leq\|r_-(z)\|_{L^\infty}\|r_+(z)\|^2_{H^{1}(\mathbb{R})}),
\end{aligned}
\end{equation}
and
\begin{equation}\label{mtea20}
\begin{aligned}
\|\langle x\rangle B_2(x)\|_{L^2(\mathbb{R}^+)}&\lesssim\|r_+(z)\|_{L^\infty(\mathbb{R})}\|\langle x\rangle{\cal{P}}^-\big({r_{-}(z)}\mathrm{e}^{2i(z+1)x}\big)\|_{L^\infty(\mathbb{R}^+,L_z^2(\mathbb{R}))}\\
&\qquad\times\|{\cal{P}}^-\big({r_{-}(z)}\mathrm{e}^{2i(z+1)x}\big)\|_{L^2(\mathbb{R}^+,L_z^2(\mathbb{R}))}\\
&\leq\|r_+(z)\|_{L^\infty}\|r_-(z)\|^2_{H^{1}(\mathbb{R})})
\end{aligned}
\end{equation}

We note, from \eqref{mte29}, that
\begin{equation}\label{mtea29}
\begin{aligned}
\|u\|_{L^2(\mathbb{R}^+)}&\leq \|A_1(x)\|_{L^2(\mathbb{R}^+)}+\|A_2(x)\|_{L^2(\mathbb{R}^+)},\\
\|\langle x\rangle u\|_{L^2(\mathbb{R}^+)}&\leq \|\langle x\rangle A_1(x)\|_{L^2(\mathbb{R}^+)}+\|\langle x\rangle A_2(x)\|_{L^2(\mathbb{R}^+)}.
\end{aligned}
\end{equation}
Using the sharp Gagliardo-Nirenberg inequality $\|u\|^3_{L^6}\leq\frac{2}{\pi}\|u\|^2_{L^2}\|u_x\|_{L^2}$, we find that
\[\begin{aligned}
\|u_x\|_{L^2(\mathbb{R}^+)}\leq& 2\|u\|_{L^2(\mathbb{R}^+)}+\frac{1}{\pi}\|u\|^2_{L^2(\mathbb{R}^+)}\|u_x\|_{L^2(\mathbb{R}^+)}\\
&\quad+\|B_1(x)\|_{L^2(\mathbb{R}^+)}+\|B_2(x)\|_{L^2(\mathbb{R}^+)}.
\end{aligned}\]
If $\|u\|_{_{L^2(\mathbb{R}^+)}}\leq\sqrt{\pi}$, we have
\begin{equation}\label{mtea30}
\|u_x\|_{L^2(\mathbb{R}^+)}\leq\frac{\pi}{\pi-\|u\|^2_{L^2(\mathbb{R}^+)}}\big(2\|u\|_{L^2(\mathbb{R}^+)}
+\|B_1(x)\|_{L^2(\mathbb{R}^+)}+\|B_2(x)\|_{L^2(\mathbb{R}^+)}\big).
\end{equation}
and
\begin{equation}\label{mtea31}
\begin{aligned}
\|\langle x\rangle u_x\|_{L^2(\mathbb{R}^+)}\leq&\frac{\pi}{\pi-\|u\|^2_{L^2(\mathbb{R}^+)}}\big(2\|\langle x\rangle u\|_{L^2(\mathbb{R}^+)}\\
&\quad+\|\langle x\rangle B_1(x)\|_{L^2(\mathbb{R}^+)}+\|\langle x\rangle B_2(x)\|_{L^2(\mathbb{R}^+)}\big).
\end{aligned}
\end{equation}
In addition,
\begin{equation}\label{mtea32}
\begin{aligned}
\|u_{xx}\|_{L^2(\mathbb{R}^+)}\leq&2\|u_x\|_{L^2(\mathbb{R}^+)}+\frac{3}{4}\|u\|^2_{H^1(\mathbb{R}^+)}\|u_x\|_{L^2(\mathbb{R}^+)}\\
&+\|\partial_xB_1(x)\|_{L^2(\mathbb{R}^+)}+\|\partial_xB_2(x)\|_{L^2(\mathbb{R}^+)}.
\end{aligned}
\end{equation}

The Lemma is proved by substituting \eqref{mtea3}-\eqref{mtea5} and \eqref{mtea33}-\eqref{mtea20} into \eqref{mtea29}-\eqref{mtea32}. \qquad $\square$

\begin{lemma}\label{lem13}
 The map
\begin{equation}
\label{b76}
{H}^{1}(\mathbb{R})\cap {L}^{2,1}(\mathbb{R})\ni(r_{-},r_{+})\longmapsto u\in {H}^{2}(\mathbb{R}^+)\cap {H}^{1,1}(\mathbb{R}^+),
\end{equation}
is Lipschitz continuous. Moreover, let $r_{\pm}, \tilde{r}_{\pm}\in {H}^{1}(\mathbb{R})\cap {L}^{2,1}(\mathbb{R})$ admitting $\| r_{\pm}\|_{{H}^{1}\cap {L}^{2,1}(\mathbb{R})}$, $\|\tilde{r}_{\pm}\|_{{H}^{1}\cap {L}^{2,1}(\mathbb{R})}\leq\rho$, the associated potentials $v,\tilde{v}$ satisfy the bound
\begin{equation}
\label{b74cd}
\| u-\widetilde{u}\|_{{H}^{2}(\mathbb{R}^{+})\cap {H}^{1,1}(\mathbb{R}^{+})}\leq C(\rho)(\| r_{+}-\tilde{r}_{+}\|_{{H}^{1}(\mathbb{R})\cap {L}^{2,1}(\mathbb{R})}+\| r_{-}-\tilde{r}_{-}\|_{{H}^{1}(\mathbb{R})\cap {L}^{2,1}(\mathbb{R})}),
\end{equation}
for some constant $C(\rho)$.
\end{lemma}
{\bf Proof}~~ The Lemma can be proved from \eqref{mte29}-\eqref{mte29b} by repeating the same estimates in Lemma \ref{lem12}. \qquad $\square$

Secondly, we consider the case $x\in\mathbb{R}^-$.
We note that, in the reconstruction formula \eqref{mte27} and \eqref{mte28},
for every $x\in\mathbb{R}^-$, $\mathrm{e}^{\pm 2i(z+1)x}$ in the representations will be difficult to control in $z\in \mathbb{C}^\pm$. So the estimate for $x\in\mathbb{R}^-$ cannot proceed.

To overcome it, we can define modified reflection coefficients
\begin{equation}\label{mte43}
r_{\pm,\delta}(z)=\overline{\delta_+(z)}\overline{\delta_-(z)}r_\pm(z),
\end{equation}
for every $z\in\mathbb{R}$, and define new jump matrix
\begin{equation}\label{mte44}
R_{\delta}(x;z)=\left(\begin{matrix}
0&\overline{r_{+,\delta}(z)}\mathrm{e}^{-2i(z+1)x}\\
r_{-,\delta}(z)\mathrm{e}^{2i(z+1)x}&\overline{r_{+,\delta}(z)}r_{-,\delta}(z)
\end{matrix}\right),
\end{equation}
where the functions
\begin{equation}\label{mte41}
\delta_{\pm}(z)=\mathrm{e}^{{\mathcal{P}}^{\pm}\log(1+\overline{r_{+}}r_{-})}, \quad z\in \mathbb{R},
\end{equation}
admit a scalar Riemann-Hilbert problem \cite{aom137-295}
\begin{equation}\label{mte39}
\begin{aligned}
&\delta_{+}(z)-\delta_{-}(z)=\overline{r_{+}(z)}r_{-}(z)\delta_{-}(z), \quad z\in \mathbb{R},\\
&\delta_{\pm}(z)\to1, \quad   z \to\infty.
\end{aligned}
\end{equation}

Now, we can define modified analytic functions
\begin{equation}\label{mte45}
 M_{\pm,\delta}(x;z)=M_\pm(x;z){\rm diag}(\delta_\pm^{-1}(z),\delta_\pm(z)),
\end{equation}
then \eqref{mtc16} can be rewritten as
\begin{equation}\label{mte46}
M_{+,\delta}(x;z)-M_{-,\delta}(x;z)=M_{-,\delta}(x;z)R_\delta(x;z), \quad z\in \mathbb{R},
\end{equation}
\begin{equation}\label{mte47}
M_{\pm,\delta}(x;z)\to I, \quad z\to\infty,
\end{equation}
The solution of \eqref{mte46} and \eqref{mte47} have the following form
\begin{equation}\label{mte48}
M_{\pm,\delta}(x;z)=I+{\cal{C}}\big(M_{-,\delta}(x;\cdot)R_\delta(x;\cdot)\big)(z), \quad z\in\mathbb{C}^\pm.
\end{equation}

Definiting $M_{\pm,\delta}(x;z)$ as
\begin{equation}\label{mte49}
 M_{\pm,\delta}(x;z)=\big(\xi_{\pm,\delta}(x;z),\eta_{\pm,\delta}(x;z)\big),
\end{equation}
from \eqref{mte46}, we have
\[\xi_{+\delta}(x;z)=\xi_{-,\delta}(x;z)+\eta_{-,\delta}(x;z)r_{-,\delta}(z)\mathrm{e}^{2i(z+1)x}.\]
Moreover, from \eqref{mte48}, we have
\begin{equation}\label{mte50}
 \lim\limits_{z\to\infty}z(\xi_{\pm,\delta}-e_1)=-\frac{1}{2\pi i}\int_{\mathbb{R}}r_{-,\delta}(z)\mathrm{e}^{2i(z+1)x}\eta_{-,\delta}(x;z)\mathrm{d}z,
\end{equation}
and
\begin{equation}\label{mte51}
 \lim\limits_{z\to\infty}z(\eta_{\pm,\delta}-e_2)=-\frac{1}{2\pi i}\int_{\mathbb{R}}\overline{r_{+,\delta}(z)}\mathrm{e}^{-2i(z+1)x}\xi_{+,\delta}(x;z)\mathrm{d}z.
\end{equation}
Then for every $x\in \mathbb{R}^-$, we have
\begin{equation}\label{mte52}
\begin{aligned}
iu&=-4\lim\limits_{z\to\infty}z\eta_-^{(1)}(x;z)=-4\lim\limits_{z\to\infty}z\eta_{-,\delta}^{(1)}(x;z)\\
&=\frac{2}{\pi i}\int_{\mathbb{R}}\overline{r_{+,\delta}(z)}\mathrm{e}^{-2i(z+1)x}\xi_{+,\delta}^{(1)}(x;z)\mathrm{d}z,
\end{aligned}
\end{equation}
\begin{equation}\label{mte53}
\begin{aligned}
i\bar{u}_{x}+2\overline{u}-\frac{1}{2}|u|^{2}\overline{u}&=-2i\lim\limits_{z\to\infty}z\xi_-^{(2)}(x;z)=-2i\lim\limits_{z\to\infty}z\xi_{-,\delta}^{(2)}(x;z)\\
&=\frac{1}{\pi }\int_{\mathbb{R}}{r_{-,\delta}(z)}\mathrm{e}^{2i(z+1)x}\eta_{-,\delta}^{(2)}(x;z)\mathrm{d}z.
\end{aligned}
\end{equation}

The limit of equation \eqref{mte48} on the real $z$ line gives to the projection equations
\begin{equation}\label{mte54}
\xi_{+,\delta}(x;z)-e_1={\cal{P}}^+\big(r_{-,\delta}(z)\mathrm{e}^{2i(z+1)x}\eta_{-,\delta}(x;z)\big),
\end{equation}
\begin{equation}\label{mte55}
\eta_{-,\delta}(x;z)-e_2={\cal{P}}^-\big(\overline{r_{+,\delta}(z)}\mathrm{e}^{-2i(z+1)x}\xi_{-,\delta}(x;z)\big),
\end{equation}
Then we have
\begin{equation}\label{mte56}
M_\delta(x;z)-{\cal{P}}^-\big(M_\delta(x;z)R_\delta(x;z)\big)=F_\delta(x;z),
\end{equation}
where
\[M_\delta(x;z)=\big(\xi_{+,\delta}(x;z)-e_1, \eta_{-,\delta}(x;z)-e_2\big)\left(\begin{matrix}
1&0\\
r_{-,\delta}(z)\mathrm{e}^{2i(z+1)x}&1
\end{matrix}\right),\]
\[F_\delta(x,z)=\big({\cal{P}}^+(r_{-,\delta}(z)\mathrm{e}^{2i(z+1)x})e_2, {\cal{P}}^-(\overline{r_{+,\delta}(z)}\mathrm{e}^{-2i(z+1)x})e_1\big).\]

For the projector operator \eqref{mtd2}, we can obtain ${\cal{P}}_++{\cal{P}}_-=-i{\cal{H}}$. Since ${\cal{H}}\log(1+\overline{r_+}r_-)$ is real according to \eqref{mtc11}, then
\begin{equation}\label{mte42}
|\delta_+(z)\delta_-(z)|=\big|\mathrm{e}^{-i{\cal{H}}\log(1+\overline{r_+}r_-)}\big|=1,
\end{equation}
So if $r_\pm(z)\in L^{2,1}(\mathbb{R})$, then $r_{\pm,\delta}(z)=\overline{\delta_+(z)\delta_-(z)}r_\pm(z)\in {L}^{2,1}(\mathbb{R})$.

We note that $r_{\pm,\delta}(z)\in {H}^{1}(\mathbb{R})$ if $r_{\pm}(z)\in {H}^{1}(\mathbb{R})$.  In fact,
\[\partial_z(\overline{\delta_-\delta_+}r_\pm)=(\overline{\delta_-\delta_+}r_\pm)
\big(-i\partial_z{\cal{H}}\log(1+\overline{r_+}r_-)\big)+\overline{\delta_-\delta_+}\partial_z(r_\pm).\]
Since $\|{\cal{H}}f\|_{L^2(\mathbb{R})}=\|f\|_{L^2(\mathbb{R})}$ and $\|\partial_z({\cal{H}}f)\|_{L^2}=\|\partial_zf\|_{L^2}$.
Thus, we find that $\overline{\delta_+(z)\delta_-(z)}r_\pm(z)\in {H}^{1}(\mathbb{R})\cap L^{2,1}(\mathbb{R})$.

Using Proposition \ref{pro2}, the statements of Lemma \ref{lem13} can be extended to the negative half-line.
\begin{lemma}\label{lem15}
If $r_{\pm}(z)\in {H}^{1}(\mathbb{R})\cap {L}^{2,1}(\mathbb{R})$ and admit the condition \eqref{mtc8}, then $u\in{H}^{2}(\mathbb{R}^-)\cap {H}^{1,1}(\mathbb{R}^-)$. In addition, we have
\begin{equation}
\label{b74bb}
\| u\|_{{H}^{2}(\mathbb{R}^{-})\cap {H}^{1,1}(\mathbb{R}^{-})}\leq C(\| r_{+,\delta}\|_{{H}^{1}(\mathbb{R})\cap {L}^{2,1}(\mathbb{R})}+\| r_{-,\delta}\|_{{H}^{1}(\mathbb{R})\cap {L}^{2,1}(\mathbb{R})},
\end{equation}
where $C$ is a positive constant depending on $\| r_{\pm,\delta}\|_{{H}^{1}(\mathbb{R})\cap {L}^{2,1}(\mathbb{R})}$.
\end{lemma}

\begin{lemma}\label{lem16}
The map
\begin{equation}
\label{b76a}
{H}^{1}(\mathbb{R})\cap {L}^{2,1}(\mathbb{R})\ni(r_{-},r_{+})\longmapsto u\in {H}^{2}(\mathbb{R}^-)\cap {H}^{1,1}(\mathbb{R}^-),
\end{equation}
is Lipschitz continuous. Moreover, let $r_{\pm}, \tilde{r}_{\pm}\in {H}^{1}(\mathbb{R})\cap {L}^{2,1}(\mathbb{R})$ admitting $\| r_{\pm}\|_{{H}^{1}\cap {L}^{2,1}(\mathbb{R})}$, $\|\tilde{r}_{\pm}\|_{{H}^{1}\cap {L}^{2,1}(\mathbb{R})}\leq\rho$, the associated potentials $v,\tilde{v}$ satisfy the bound
\begin{equation}
\label{b74cc}
\| u-\widetilde{u}\|_{{H}^{2}(\mathbb{R}^{-})\cap {H}^{1,1}(\mathbb{R}^{-})}\leq C(\rho)(\| r_{+}-\tilde{r}_{+}\|_{{H}^{1}(\mathbb{R})\cap {L}^{2,2}(\mathbb{R})}+\| r_{-}-\tilde{r}_{-}\|_{{H}^{1}(\mathbb{R})\cap {L}^{2,2}(\mathbb{R})}),
\end{equation}
for some constant $C(\rho)$.
\end{lemma}
\setcounter{equation}{0}
\section{Global solutions to the MTNLS-GI equation}\label{sec7}

In the following, we firstly consider the time evolution of the scattering data.
Define the fundamental solution for every $t\in[0,T]$
\[\begin{aligned}
&\widehat{Y}(x,t;k)=\mathrm{e}^{i\theta(x,t;k)}\varphi_{\pm}(x,t;k),\\
&\widehat{Y}(x,t;k)=\mathrm{e}^{-i\theta(x,t;k)}\phi_{\pm}(x,t;k),
\end{aligned}\]
with
\[
\varphi_{\pm}(x,t;k)\to e_{1}, \quad \phi_{\pm}(x,t;k)\to e_{2},\quad x\to\pm\infty.
\]
Here
\[\theta(x,t;k)=(k^{2}+1) x-2(k^{2}+1)^{2}t\]
From \eqref{mtb44},\eqref{mtb44a}, the bounded function $\varphi_{\pm},\phi_{\pm}$ satisfy similar analytical properties, and satisfy the scattering relationship
\begin{equation}
\label{c1}
\varphi_{-}(t,x;k)=a(k)\varphi_{+}(t,x;k)+b(k)\mathrm{e}^{2i\theta(x,t;k)}\phi_{+}(x,t;k),\quad x\in R,k\in R\cup iR.
\end{equation}
Furthermore, the time-dependent reflection coefficients can be obtained
\begin{equation}
\label{c3}
r_{\pm}(t;z)=r_{\pm}(z)\mathrm{e}^{2i(z+1)^{2}t}.
\end{equation}
where $r_{\pm}(0,k)$ is the initial spectral data related to initial condition $u(0,x)$.
Lemma \ref{lem12}-Lemma \ref{lem16} can be extended to the reflection coefficients $r_{\pm}(t;z)$ by similar discussions.

\begin{proposition}\label{pro7}
For every $u_0\in H^2(\mathbb{R})\cap H^{1,1}(\mathbb{R})$ such that the linear spectral problem \eqref{a2} admits no eigenvalues or resonances. For every $t\in[0,T]$, there exists a unique local solution
\[\begin{array}{c}
u(x,t)\in C\big(H_x^2(\mathbb{R})\cap H_x^{1,1}(\mathbb{R}), [0,T]\big),\\
\end{array}
\]
to the Cauchy problem of NLS-GI equation \eqref{a1}. Moreover, the map
\begin{equation}\label{mtf1}
\begin{array}{c}
H^2(\mathbb{R})\cap H^{1,1}(\mathbb{R})\ni u_0\\
\end{array}
\mapsto \begin{array}{c}
u(x,t)\in C\big(H_x^2(\mathbb{R})\cap H_x^{1,1}(\mathbb{R}), [0,T]\big),\\
\end{array}
\end{equation}
is Lipschitz continuous.
\end{proposition}
{\bf Proof}.~~We note that, the potential $u(x,t)$ can be controlled by $\|r_\pm(t;z)\|_{H_z^1(\mathbb{R})\cap L_z^{2,1}(\mathbb{R})}$ via the same lemmas as Lemma \ref{lem12} and Lemma \ref{lem15}.  Furthermore, from \eqref{c3}, we find that
$\|r_\pm(t;z)\|_{L^2(\mathbb{R})}=\|r_\pm(z)\|_{L^2(\mathbb{R})}, \ \|r_\pm(t;z)\|_{L^{2,1}(\mathbb{R})}=\|r_\pm(z)\|_{L^{2,1}(\mathbb{R})},$
and
\[\|\partial_zr_\pm(t;z)\|_{L^2(\mathbb{R})}\leq\|\partial_zr_\pm(z)\|_{L^2(\mathbb{R})}+4T\|r_\pm(z)\|_{L^2(\mathbb{R})}+4T\|r_\pm(z)\|_{L^{2,1}(\mathbb{R})}.\]
Thus, the potential $u(x,t)$ can be further controlled by $C(T)\|r_\pm(z)\|_{H_z^1(\mathbb{R})\cap L_z^{2,1}(\mathbb{R})}$ in view of the Lemma \ref{lem8}.
As a result, $u(x,t)$ can be covered by $c(T)\|u_0(x)\|_{H^2(\mathbb{R})\cap H^{1,1}(\mathbb{R})}$, where the positive constant $c(T)$ may grow at most polynomially in $T$ but is remains finite for every $T>0$. In fact,
\[\begin{aligned}
\|u(x,t+\Delta t)\|_{H^2(\mathbb{R})\cap H^{1,1}(\mathbb{R})}&\leq C\|r_\pm(z)(\mathrm{e}^{2i(z+1)^{2}(t+\Delta t}-\mathrm{e}^{2i(z+1)^{2}t})\|_{H_z^1(\mathbb{R})\cap L_z^{2,1}(\mathbb{R})}\\
&\leq C|\Delta t|\|r_\pm(z)\|_{H_z^1(\mathbb{R})\cap L_z^{2,1}(\mathbb{R})}.
\end{aligned}\]

According to discussions in the previous sections, we get that there exists a unique local solution
\[\begin{array}{c}
u(x,t)\in C\big(H_x^2(\mathbb{R})\cap H_x^{1,1}(\mathbb{R}), [0,T]\big),\\
\end{array}
\]
to the Cauchy problem of the NLS-GI equation \eqref{a1}.  \qquad $\square$

\begin{theorem}
For every $u_0\in H^2(\mathbb{R})\cap H^{1,1}(\mathbb{R})$ such that the linear spectral problem \eqref{a2} admits no eigenvalues or resonances. For every $t\in[0,\infty)$, there exists a unique global solution
\[\begin{array}{c}
u(x,t)\in C\big(H_x^2(\mathbb{R})\cap H_x^{1,1}(\mathbb{R}), [0,\infty)\big),\\
\end{array}
\]
to the Cauchy problem of NLS-GI equation \eqref{a1}. Moreover, the map
\begin{equation}\label{mtf2}
\begin{array}{c}
H^2(\mathbb{R})\cap H^{1,1}(\mathbb{R})\ni u_0\\
\end{array}
\mapsto \begin{array}{c}
u(x,t)\in C\big(H_x^2(\mathbb{R})\cap H_x^{1,1}(\mathbb{R}), [0,\infty)\big),\\
\end{array}
\end{equation}
is Lipschitz continuous.
\end{theorem}
{\bf Proof}~~ Assuming the existencen of a constant $T_{max}<\infty$
\[\lim\limits_{t\to T_{m}}\big(\|u(x,t)\|_{H_x^2\cap H_x^{1,1}}\big)=\infty,\]
On the other hand, from the proof of the Proposition \ref{pro7}, $\|u(x,t)\|_{H_x^2\cap H_x^{1,1}}$ is control by $\|u_0(x)\|_{H_x^2\cap H_x^{1,1}}$ which is bound for any $t\in[0,T]$. This contradicts the assumption. Thus the local solution $u$ can be continued globally in time for every $T>0$.   \qquad $\square$

\section*{Statements and Declarations}

\noindent{\bf Conflict of interest} The authors declare that they have no conflicts of interest.

\noindent{\bf Funding} This work is supported by the National Natural Science
Foundation of China under Grant (No. 12571268).


\end{document}